\documentclass[a4paper,12pt]{article}

\usepackage[english]{babel}
\usepackage{amsfonts,amsmath,amsthm,mathrsfs,amscd}
\usepackage{dsfont}
\usepackage{cite}
\usepackage{authblk}

\usepackage{enumitem}
\setlist{nolistsep, itemsep=0.1cm,parsep=0pt}

\usepackage{graphicx}

\newtheorem{theorem}{Theorem}
\newtheorem{lemma}{Lemma}

\newtheorem{prop}{Proposition}

\theoremstyle{remark}
\newtheorem{remark}{Remark}

\providecommand{\keywords}[1]
{
  \small	
  \textit{Key Words:~} #1
}

\providecommand{\MSC}[1]
{
  \small	
  \textit{Mathematics Subject Classification 2010:~} #1
}

\begin{document}

\title{Description of random fields\\by systems of conditional distributions}

\author{Khachatryan L.A.}

\affil{\small Institute of Mathematics, NAS of RA, {\it linda@instmath.sci.am}}

\date{}

\maketitle

\begin{abstract}
In this paper, we consider the direct and inverse problems of the description of lattice positive random fields by various systems of finite-dimensional (as well as one-point) probability distributions parameterized by boundary conditions. In the majority of cases, we provide necessary and sufficient conditions for the system to be a conditional distribution of a (unique) random field. The exception is Dobrushin-type systems for which only sufficient conditions are known. Also, we discuss possible applications of the considered systems.\end{abstract}

\keywords{Random field, conditional distribution, specification, Markov random field}

\MSC{60G60, 60E05, 60J99}

\section*{Introduction}

A random field is a probability measure on the infinite-dimensional space of its realizations. This definition, which is of great theoretical importance, is not always suitable for applications. From a practical point of view, it is much more convenient to deal with systems of finite-dimensional (conditional or unconditional) distributions equivalently describing a random field.

Systems of finite-dimensional distributions can be subdivided into two types: systems generated by a random field and those given autonomously through the main system-forming property of its elements --- consistency conditions. Accordingly, two problems come to the fore. The direct problem is the problem of unique determination (restoration) of a random field by the system generated by it. The inverse problem is the problem of the existence of a random field with an \textit{a~priori} given system of distributions.

Kolmogorov~\cite{Kolm} was the first one who considered both the direct and the inverse problems of the description of random processes by a system of unconditional finite-dimensional distributions. Over time, it became clear that in many cases, it is convenient to impose restrictions not on unconditional distributions but on their relations, that is, on conditional distributions.

The idea of specifying a random object through conditional probabilities is very old and goes back to the concept of a Markov chain. At the end of the 1960s, the approach to describing random objects by systems of conditional distributions was further developed. In connection with mathematical problems of statistical physics, Dobrushin~\cite{D2} and, independently, Lanford and Ruelle~\cite{LR} introduced the fundamental concept of a Gibbs random field. Its definition is based on the idea of specifying a random field through a special system of finite-dimensional distributions parameterized by boundary conditions --- Gibbs specification.

Dobrushin's approach was further developed by Dachian and Nahapetian in~\cite{DN1998, DN2001, DN2004}, where it was shown that the description of Dobrushin-type specifications (and, consequently, random fields) can be carried out by systems of consistent one-point distributions parameterized by infinite boundary conditions. Systems of one-point distributions on general measurable spaces were studied by Fernandez and Maillard~\cite{FernMill}. The problem of describing a random field by a system of consistent one-point distributions with finite boundary conditions was solved by Dalalyan and Nahapetian~\cite{DN2011}.

In this paper, we consider various systems of conditional probability distributions, study their properties, and present for such systems solutions for both direct and inverse problems. Some of the systems are fully studied for the first time: the system of finite-dimensional distributions parameterized by finite boundary conditions (Section~\ref{Section-Qf}), the system of Palm-type distributions (Section~\ref{Section-Q-Pi}), and the systems of finite-dimensional and one-point distributions parameterized by various boundary conditions (Section~\ref{Section-Q}). For the system of one-point distributions with finite boundary conditions, we slightly improve known results (Section~\ref{Section-Qf-1}). Also, we consider the systems of finite-dimensional and one-point distributions parameterized by Dobrushin's type boundary conditions and formulate the corresponding results within the framework of the approach developed in the paper  (Section~\ref{Section-QD}).

The work is mainly theoretical, but possible areas of application of the obtained results will be outlined. Some of the statements of this paper (with proof ideas) are given in~\cite{KhN20} (see also~\cite{AN}). Here, for the sake of simplicity, we restrict ourselves to considering only positive random fields with a finite phase space. However, the results can be carried over more general settings.

\section{Preliminaries}

Let $X \subset \mathbb{R}$ be a non-empty finite set, $1 < \vert X \vert < \infty$, and let $\mathbb{Z}^d$ be a $d$-dimensional integer lattice (the set of $d$-dimensional vectors with integer components), $d \geq 1$.

For $S \subset \mathbb{Z}^d$, denote by $W(S) = \{ V \subset S: \left| V \right| < \infty \}$ the set of all finite subsets of $S$, and let $W_n(S) =  \{V \subset S: \vert V \vert = n\}$ be the set of $n$-element subsets of $S$, $n \ge 1$. When $S = \mathbb{Z}^d$, we will use simpler notations $W$ and $W_n$, respectively. For convenience of notations, we will usually omit braces for one-point sets $\{t\}$, $t \in \mathbb{Z}^d$.

\enlargethispage{\baselineskip}

A neighborhood system in $\mathbb{Z}^d$ is a system $\partial = \{ \partial t, t \in \mathbb{Z}^d \}$ of finite-dimensional subsets $\partial t$ of the lattice $\mathbb{Z}^d$ such that $t \notin \partial t$ and $s \in \partial t$ if and only if $t \in \partial s$, $s \in \mathbb{Z}^d$. If a neighborhood system $\partial$ is defined, for any $V \in W$, put $\partial V = \{s \in \mathbb{Z}^d \backslash V: s \in \partial t, t \in V\}$.

Denote by $X^S = \{ x = ( x_t ,t \in S ): x_t \in X^{\{t\}}, t \in S \}$ the set of configurations on $S$,  $S \subset \mathbb{Z}^d$, i.e., the set of functions defined on $S$ with values in $X$. If $S = {\O}$, we assume that $X^{\O} = \{ \boldsymbol{{\O}} \}$ where $\boldsymbol{{\O}}$ is an empty configuration. For any $S, T \subset \mathbb{Z}^d$ such that $S \cap T = {\O}$ and any $x \in X^S$ and $y \in X^T$, denote by $xy$ the concatenation of $x$ and $y$, that is, the configuration on $S \cup T$ coinciding with $x$ on $S$ and with $y$ on $T$. We assume that concatenation of $x$ with an empty configuration $\boldsymbol{{\O}}$ coincides with $x$, i.e., $x \boldsymbol{{\O}} = x$ for all $x \in X^S$, $S \subset \mathbb{Z}^d$. Finally, for $T \subset S$, denote by $x_T$ the restriction of the configuration $x \in X^S$ to $T$.

When some enumeration $V=\{t_1, t_2, ..., t_n\}$ of the points of $V \in W$ is fixed, for brevity, we will denote by $(xu)_j$ the concatenation of configurations $x_{\{t_1,...,t_{j-1}\}}$ and $u_{\{t_{j+1,...,t_n}\}}$, that is
\begin{equation}
	\label{xu-notations}
	\begin{array}{c}
		(xu)_j = x_1 ... x_{j-1} u_{j+1} ... u_n, \quad 1 < j < n, \\
		\\
		(xu)_1 = u_2 u_3 ... u_{t_n}, \qquad (xu)_n = x_1 x_2 ... x_{n-1},
	\end{array}
\end{equation}
where $x_j = x_{t_j}$, $u_j = u_{t_j}$, $1 \le j \le n$, $n = \vert V \vert$, and $x,u \in X^V$.

For $x \in X^S$, $S \subset \mathbb{Z}^d$, we call the set $S$ the support of configuration $x$ and denote it by $s(x)$. For any $V \in W$, for sets of configurations with supports not intersecting with $V$ (or, simply, for configurations outside $V$), we will use the following notations
$$
\widehat{X_V} = \bigcup \limits_{{\O} \neq S \subset \mathbb{Z}^d \backslash V} X^S, \qquad
\widehat{ X_V^f} = \bigcup \limits_{{\O} \neq S \in W(\mathbb{Z}^d \backslash V)} X^S, \qquad
\widehat{ X_V^{f0}} = \bigcup \limits_{S \in W(\mathbb{Z}^d \backslash V)} X^S.
$$
It is clear that $\widehat{X_V^f} \subset \widehat{X_V}$ and $\widehat{X_V^f} \subset \widehat{X_V^{f0}}$.

For a sequence $\Lambda = \{ \Lambda_n \}_{n \ge 1}$ of sets $\Lambda_n \in W$ and for $S \subset \mathbb{Z}^d$, the notation $\Lambda_n \uparrow S$ means that the sequence $\Lambda$ is increasing and converges to $S$, i.e., ${\Lambda_n \subset \Lambda_{n+1}}$ and $\bigcup \limits_{n=1}^\infty \Lambda_n = S$. For a family $\{g_\Lambda, \Lambda \in W(S)\}$ of functions, the notation $\lim \limits_{\Lambda \uparrow S} g_\Lambda(x_\Lambda) = a$, $x \in X^S$, means that for any increasing sequence $\{\Lambda_n\}_{n \ge 1}$ of finite sets converging to $S$, we have $\lim \limits_{n \to \infty} g_{\Lambda_n}(x_{\Lambda_n}) = a$. A real-valued function $g$ on $X^S$ is called quasilocal if
$$
\lim_{\Lambda \uparrow S} \sup_{x,y \in X^S \,:\, x_\Lambda = y_\Lambda} \bigl| g(x) - g(y) \bigr| = 0.
$$

For any $S \subset \mathbb{Z}^d$, denote by $\mathscr{B}^S$ the $\sigma$-algebra generated by cylinder subsets of $X^S$. In the case of a finite subset of the lattice $V \in W$, $\mathscr{B}^V$ is the $\sigma$-algebra of all subsets of $X^V$.

Probability distribution on $(X^S, \mathscr{B}^S)$ will be denoted by Latin letters indexed by the set $S$, for example, $P_S$. In the case $S = \mathbb{Z}^d$, the subscript will be omitted, i.e., we will write $P$ instead of $P_{\mathbb{Z}^d}$. If $S = {\O}$, there exists only one probability distribution $P_{{\O}} (\boldsymbol{{\O}}) = 1$. For $T \subset S \subset \mathbb{Z}^d$ and any probability distribution $P_S$ on $(X^S, \mathscr{B}^S)$, denote by $P_T$ its restriction $(P_S)_T$ on $(X^T, \mathscr{B}^T)$. In the case of finite subsets of the lattice $V \subset \Lambda \in W$, one has
$$
P_V (x) = (P_\Lambda)_V (x) = \sum \limits_{y \in X^{\Lambda \backslash V}} P_\Lambda(xy), \qquad x \in X^V.
$$

A probability distribution $P$ on $(X^{\mathbb{Z}^d}, \mathscr{B}^{\mathbb{Z}^d})$ is called a random field. For a random field $P$, the set of probability distributions $\{ P_V ,V \in W \}$ with $P_V = (P)_V$ is called its system of finite-dimensional (unconditional) distributions.

It is well known (see, for instance, the classical work~\cite{Kolm} by Kolmogorov) that any random field is restored by its system of finite-dimensional unconditional distributions. In this regard, we will often identify a random filed $P$ with the system of its finite-dimensional distributions  and will write $P = \{P_V, V \in W\}$. An autonomously given system of finite-dimensional distributions $\{{\rm P}_V, V \in W\}$ satisfying Kolmogorov's consistency condition $({\rm P}_\Lambda)_V = {\rm P}_V$ for all $V \subset \Lambda \in W$ uniquely determines a random field $P$ such that $(P)_V = {\rm P}_V$, $V \in W$.

A random field $P$ is called positive if its system of finite-dimensional distributions is positive, that is, for any $V \in W$, one has $P_V (x) > 0$ for all $x \in X^V$. In the framework of this paper, we will consider only positive random fields.

For a random field $P$, a conditional probability $Q_V^z$ on $X^V$ under a finite boundary condition $z \in X^S$, $S \in W(\mathbb{Z}^d \backslash V)$, is defined as follows:
\begin{equation}
	\label{P->QP}
	Q_V^z (x) = \dfrac{P_{V \cup S}(xz)}{P_S(z)}, \qquad x \in X^V, \, V \in W.
\end{equation}
In the case of infinite boundary condition $z \in X^S$, $S \subset \mathbb{Z}^d \backslash V$, put
\begin{equation}
	\label{P->QP-lim}
	Q_V^z (x) = \lim \limits_{\Lambda \uparrow S} Q_V^{z_\Lambda}(x) = \lim \limits_{\Lambda \uparrow S} \dfrac{P_{V \cup \Lambda}(xz_\Lambda)}{P_\Lambda(z_\Lambda)}, \qquad x \in X^V, V \in W,
\end{equation}
where the limit exists for almost all (with respect to the measure $P$) configurations $z$. The term ``boundary condition'' for the configuration $z$ is used in mathematical statistical physics. In the future, we will adhere to this terminology.

Finally, note that conditional probabilities of a random field satisfy the following Sullivan's inequalities (see~\cite{Sulliv})
\begin{equation}
	\label{Q-ineq}
	\inf \limits_{y \in X^S: y_\Lambda = z} Q_V^y (x) \le Q_V^z (x) \le \sup \limits_{y \in X^S: y_\Lambda = z} Q_V^y (x),
\end{equation}
where $V, \Lambda \in W$, $\Lambda \subset S \subset \mathbb{Z}^d \backslash V$ and $x \in X^V$, $z \in X^\Lambda$.

\section{Systems of conditional distributions. The direct and inverse problems}

Further, we will study various systems of finite-dimensional probability distributions parame\-te\-rized by boundary conditions. Their general structure (regardless of whether they were generated by a random field or not) has the following form:
$$
{\rm Q} = \{q_V^z, z \in \Upsilon_V, V \in \mathscr{W}\},
$$
where $\Upsilon_V \subset \widehat{X_V}$ is the set of admissible boundary conditions defined outside~$V$, $\mathscr{W} \subset W$ is a family of finite subsets of the lattice $\mathbb{Z}^d$, and for any fixed $V$ and boundary condition $z \in \Upsilon_V$, the function $q_V^z$ is a probability distribution on~$X^V$.

Any such system of probability distributions is specified by two sets: $\mathscr{W}$ and $\Upsilon = \{ \Upsilon_V, V \in \mathscr{W}\}$. As $\mathscr{W}$, one can consider, for example, the set $W$ of all finite subsets of the lattice $\mathbb{Z}^d$ or the set $W_n$ of its $n$-element subsets ($n \ge 1$). A special place here is occupied by the case $n=1$, the system of sets $W_1 = \{\{t\}, t \in \mathbb{Z}^d\}$ corresponding to which is the collection of all lattice nodes. As boundary conditions, one can consider the set $\Upsilon^f = \{ \widehat{X_V^f}, V \in \mathscr{W}\}$ of configurations with finite supports or the set $\Upsilon = \{ \widehat{X_V}, V \in \mathscr{W}\}$ of configurations admitting infinite supports.

The main systems considered in the present paper are the system ${\rm Q}^f = \{q_V^z, z \in \widehat{X_V^f}, V \in W\}$ of finite-dimensional distributions parameterized by finite boundary conditions and the system ${\rm Q} = \{q_V^z, z \in \widehat{X_V}, V \in W\}$ which is the completion of the system ${\rm Q}^f$ by distributions with infinite boundary conditions.

All the other systems studied are the subsystems of the mentioned systems. For example, we will consider a Palm-type system ${\rm Q}^{\Pi} = \{q_V^z, z \in X^{\{t\}}, t\in \mathbb{Z}^d \backslash V, V \in W\}$ and Dobrushin-type system ${\rm Q}^D = \{q_V^z, z \in X^{\mathbb{Z}^d \backslash V},\linebreak {V \in W}\}$. Also, we will consider one-point systems ${\rm Q}_1^f = \{q_t^z, z \in \widehat{X^f_t}, t \in \mathbb{Z}^d\}$, ${\rm Q}_1 = \{q_t^z, z \in \widehat{X_t}, t \in \mathbb{Z}^d\}$ and ${\rm Q}_1^D = \{q_t^z, z \in X^{\mathbb{Z}^d \backslash \{t\}}, t \in \mathbb{Z}^d\}$.

The relationships of the studied systems are shown in the following diagram.
\begin{center}
	\begin{tabular}{ccccccccc}
		${\rm Q}^{\Pi}$ & $\subset$ & ${\rm Q}^f$ & $\subset$ & ${\rm Q}$ & $\supset$ & ${\rm Q}^D$ \\
		& & $\cup$ & & $\cup$ & & $\cup$ \\
		& & ${\rm Q}_1^f$ & $\subset$ & ${\rm Q}_1$ & $\supset$ & ${\rm Q}_1^D$
	\end{tabular}
\end{center}

Systems generated by a random field $P$ will be denoted by $Q_P$ or $Q(P)$. The natural requirement for such systems is that $Q_P$ restores $P$. In the case it is necessary to emphasize that the random field $P$ is restored by $Q_P$, we will use the notation $P_{Q_P}$. For a given random field $P$, we call the problem of the existence of a system $Q_P$ such that $P_{Q_P} = P$ a \emph{direct problem} of the description of random fields. A system $Q_P$ is a solution to the direct problem for the given random field $P$ if $P_{Q_P} = P$. Note that for any random field, there may exist many (equivalent) solutions to the direct problem.

For an \textit{a priori} given system ${\rm Q}$ of finite-dimensional probability distributions, we will call by an \emph{inverse problem} of the description of random fields the problem of the existence of a random field $P$ such that $Q_P = {\rm Q}$. A random field $P$ is a solution to the inverse problem for a given system ${\rm Q}$ if $Q_P = {\rm Q}$. For the system ${\rm Q}$, any solution to the inverse problem will be denoted by $P_{{\rm Q}}$. In this case, ${Q}_{P_{\rm Q}} = {\rm Q}$, and we will say that ${\rm Q}$ defines the random field~$P_{\rm Q}$. If $P_{\rm Q}$ is unique, then we will say that ${\rm Q}$ specifies it. A random field $P$ for which $Q_P = {\rm Q}$ will be called \emph{compatible} with the system~${\rm Q}$.

The solution of the direct problem makes it possible to define various classes of random fields by imposing corresponding restrictions. For example, under suitable conditions, Kolmogorov's system defines classes of Gaussian random fields, processes with independent increments or stationary random processes, while the restrictions on the systems of conditional probabilities lead to such important classes of random fields as Markov and Gibbs random fields, martingales, etc. The solution of the inverse problem provides the possibility to construct models of random fields with required properties.

Finally, note that both direct and inverse problems can be formulated not only for the pair ``random field'' -- ``system of finite-dimensional distributions'', but also for the pair ``system'' -- ``subsystem''.

\section{Systems of distributions with finite boundary conditions}

We start by considering systems of finite-dimensional distributions parameterized by finite boundary conditions. We will show that each of such systems specifies (uniquely determines) compatible with it random field. The general case (Subsection~\ref{Section-Qf}) and the Palm-type distributions (Subsection~\ref{Section-Q-Pi}) are studied in full in the present paper for the first time. One-point distributions with finite boundary conditions (Subsection~\ref{Section-Qf-1}) were the subject of the work~\cite{DN2011} by Dalalyan and Nahapetian.


\subsection{System of finite-dimensional distributions with finite boundary conditions}
\label{Section-Qf}

Let $P = \{P_V, V \in W\}$ be a random field and let $Q_P^f = \{Q_V^z, z \in \widehat{X_V^f}, V \in W\}$ be a system of its conditional probabilities (distributions) with finite boundary conditions (see~\eqref{P->QP}). The system $Q_P^f$ will be called the \emph{finite-conditional distribution of the random field $P$} or, in short, its \emph{$f$--distribution}.

From the probabilistic point of view, the definition of the system $Q_P^f$ is quite natural. This system is mentioned, for example, in~\cite{DN2009}.

Let us show that $Q_P^f$ is a solution to the direct problem for the random field $P$.

\begin{theorem}
	\label{Th-QPfin->P}
	Any random field is restored by its $f$--distribution.
\end{theorem}
\begin{proof}
	It is enough to note that finite-dimensional distributions $\{P_V, V \in W\}$ of the random field  $P$ and its finite-conditional probabilities $Q_P^f = \{Q_V^z, z \in \widehat{X_V^f}, V \in W\}$ are connected by the following relation
	\begin{equation}
		\label{QP->P}
		P_V (x) = \left(\sum \limits_{y \in X^I} \frac{Q_I^x (y)}{Q_V^y (x)} \right)^{-1}, \qquad x \in X^V, V \in W,
	\end{equation}
	where $I \in W(\mathbb{Z}^d \backslash V)$. Indeed, taking~\eqref{P->QP} into account, we can write
	$$
	\sum \limits_{y \in X^I} \frac{Q_I^x (y)}{Q_V^y (x)} = \sum \limits_{y \in X^I} \frac{P_I (y)}{P_V(x)} = \frac{1}{P_V(x)}.
	$$
\end{proof}


To solve the inverse problem associated with the system $Q_P^f$, first of all, it is necessary to answer the following question: does the system $Q_P^f$ possess such properties (consistency conditions) which allow restoring the random field $P$ without taking into account the fact that the elements of $Q_P^f$ are generated by $P$? If such characterizing properties are found, one can expect that for an \textit{a~priori} given system ${\rm Q}^f$ of distributions satisfying these properties, the inverse problem will have a solution. That is, there will exist a random field $P_{Q^f}$, the $f$--distribution $Q^f (P_{Q^f})$ of which coincides with ${\rm Q}^f$.

As it will be shown below, for the system $Q_P^f$, such characterizing property is the following one: for any disjoint sets $V,I \in W$ and boundary conditions $z \in \widehat{X_{V\cup I}^f}$, it holds
\begin{equation}
	\label{cond-prop-main}
	Q_{V \cup I}^z (xy) = Q_V^z(x) Q_I^{z x} (y),  \qquad x \in X^V, \, y \in X^I.
\end{equation}
The verification of these relations for the conditional probabilities of a random field is trivial.

A system ${\rm Q}^f = \{q_V^z, z \in \widehat{X_V^f}, V \in W\}$ of strictly positive probability distributions $q_V^z$ on $X^V$ parameterized by finite boundary conditions $z\in \widehat{X_V^f}$, $V \in W$, will be called a \emph{specification with finite boundary conditions} (or \emph{$f$--specification}) if its elements satisfy the following consistency conditions: for all disjoint sets $V,I \in W$ and all configurations $x \in X^V$, $y \in X^I$, $z \in \widehat{X_{V \cup I}}$, it holds
\begin{equation}
	\label{Qf-sogl}
	q_{V \cup I}^z (xy) = q_V^z(x) q_I^{zx} (y).
\end{equation}

It is not difficult to see that $f$--distribution $Q_P^f$ of a random field $P$ forms an $f$--specification.

Note that the positivity condition imposed on the elements of the considered system is due to the fact that in this paper, we concentrate on the problem of the description of positive random fields. The inverse problem can be solved under less restrictive positivity conditions using the ideas applied in~\cite{DN2004} (see also~\cite{KhN22}).

The following result takes place.

\begin{theorem}
	\label{Th-Qfin->P}
	Any $f$--specification specifies compatible with it random field.
\end{theorem}

To prove this theorem, we need the following properties of the elements of $f$--specification.

\begin{lemma}
	Let ${\rm Q}^f = \{ q_V^z, z \in \widehat{X_V^f}, V \in W \}$ be an $f$--specification. Then for any disjoint sets $V,I \in W$ and all configurations $x,u \in X^V$, $y,v \in X^I$, $z \in \widehat{X_{V \cup I}^{f0}}$, it holds
	\begin{equation}
		\label{Qf-4x4}
		q_V^{zy} (x) q_I^{zx} (v) q_V^{zv} (u) q_I^{zu} (y) = q_V^{zy} (u) q_I^{zu} (v) q_V^{zv} (x) q_I^{zx} (y).
	\end{equation}
	Further, for any pairwise disjoint sets $V,I, J \in W$ and any configurations $x \in X^V$, $y \in X^I$, $w \in X^J$, $z \in \widehat{X_{V \cup I \cup J}^{f0}}$, one has
	\begin{equation}
		\label{Qf-3x3}
		q_V^{zw}(x) q_I^{zx}(y) q_J^{zy}(w) = q_I^{zw}(y) q_V^{zy}(x) q_J^{zx}(w) .
	\end{equation}
\end{lemma}
\begin{proof}
	Tacking into account the consistency of the elements of ${\rm Q}^f$, for any disjoint sets $V,I \in W$ and all configurations $x,u \in X^V$, $y,v \in X^I$, $z \in \widehat{X_{V \cup I}^f}$, we can write
	$$
	\begin{array}{l}
		q_V^{zy} (x) q_I^{zx} (v) q_V^{zv} (u) q_I^{zu} (y) = \dfrac{q_{V \cup I}^z(xy)}{q_I^z(y)} \cdot \dfrac{q_{V \cup I}^z(xv)}{q_V^z(x)} \cdot \dfrac{q_{V \cup I}^z(uv)}{q_I^z(v)} \cdot \dfrac{q_{V \cup I}^z(uy)}{q_V^z(u)} =\\
		\\
		= \dfrac{q_{V \cup I}^z(uy)}{q_I^z(y)} \cdot \dfrac{q_{V \cup I}^z(uv)}{q_V^z(u)} \cdot \dfrac{q_{V \cup I}^z(xv)}{q_I^z(v)} \cdot \dfrac{q_{V \cup I}^z(xy)}{q_V^z(x)}= q_V^{zy} (u) q_I^{zu} (v) q_V^{zv} (x) q_I^{zx} (y).
	\end{array}
	$$
	To verify relation~\eqref{Qf-4x4} in the case $z = \boldsymbol{\O}$, note that due to~\eqref{Qf-sogl}, for any pairwise disjoint sets $V,I, \Lambda \in W$ and any configurations $x \in X^V$, $y \in X^I$, $z \in X^\Lambda$, we have
	\begin{equation}
		\label{Qf-fraction}
		\frac{q_{V \cup \Lambda}^y(xz)}{q_{I \cup \Lambda}^x(yz)} = \frac{q_V^y(x)}{q_I^x(y)}.
	\end{equation}
	Then, using~\eqref{Qf-fraction},~\eqref{Qf-sogl} and~\eqref{Qf-4x4} with $z \neq \boldsymbol{\O}$, we can write
	$$
	\begin{array}{l}
		\dfrac{q_V^y(x)}{q_I^x(y)} \cdot \dfrac{q_V^v(u)}{q_I^u(v)} = \dfrac{q_{V \cup \Lambda}^y(xz)}{q_{I \cup \Lambda}^x(yz)} \cdot \dfrac{q_{V \cup \Lambda}^v(uz)}{q_{I \cup \Lambda}^u(vz)} =  \dfrac{q_\Lambda^y(z) q_V^{zy}(x)}{q_\Lambda^x(z) q_I^{zx}(y)} \cdot \dfrac{q_\Lambda^v(z) q_V^{zv}(u)}{q_\Lambda^u(z) q_I^{zu}(v)} = \\
		\\
		= \dfrac{q_\Lambda^y(z) q_V^{zy}(u)}{q_\Lambda^u(z) q_I^{zu}(y)} \cdot \dfrac{q_\Lambda^v(z) q_V^{zv}(x)}{q_\Lambda^x(z) q_I^{zx}(v)} = \dfrac{q_{V \cup \Lambda}^y(uz)}{q_{I \cup \Lambda}^u(yz)} \cdot \dfrac{q_{V \cup \Lambda}^v(xz)}{q_{I \cup \Lambda}^x(vz)} = \dfrac{q_V^y(u)}{q_I^u(y)} \cdot \dfrac{q_V^v(x)}{q_I^x(v)}.
	\end{array}
	$$
	
	Further, from the consistency conditions~\eqref{Qf-sogl}, it follows that for any pairwise disjoint sets $V,I, J \in W$ and any configurations $x \in X^V$, $y \in X^I$, $w \in X^J$, $z \in \widehat{X_{V \cup I \cup J}^f}$, we have
	$$
	\begin{array}{l}
		q_V^{zw}(x) q_I^{zx}(y) q_J^{zy}(w) = \dfrac{q_{V \cup J}^z(xw)}{q_J^z(w)} \cdot \dfrac{q_{V \cup I}^z(xy)}{q_V^z(x)} \cdot \dfrac{q_{I \cup J}^z(yw)}{q_I^z(y)} = \\
		\\
		= \dfrac{q_{I \cup J}^z(yw)}{q_J^z(w)} \cdot \dfrac{q_{V \cup I}^z(xy)}{q_I^z(y)} \cdot \dfrac{q_{V \cup J}^z(xw)}{q_V^z(x)} = q_I^{zw}(y) q_V^{zy}(x) q_J^{zx}(w).
	\end{array}
	$$
	To verify relation~\eqref{Qf-3x3} in the case $z = \boldsymbol{\O}$, note that due to~\eqref{Qf-sogl}, it follows that for any disjoint sets $V,I \in W$ and all configurations $x \in X^V$, $y \in X^I$, $z \in \widehat{X_{V \cup I}}$, it holds
	\begin{equation}
		\label{Q-sogl-ekv1}
		q_V^z(x) q_I^{zx} (y) = q_I^z(y) q_V^{zy} (x).
	\end{equation}
	Then for any pairwise disjoint sets $V,I, J \in W$ and any configurations $x \in X^V$, $y \in X^I$, $w \in X^J$, we can write
	$$
	\begin{array}{c}
		q_V^w(x) q_I^{wx}(y) = q_I^w(y) q_V^{wy}(x),\\
		q_J^y(w) q_V^{yw}(x) = q_V^y(x) q_J^{yx}(w),\\
		q_I^x(y) q_J^{xy}(w) = q_J^x(w) q_I^{xw}(y).
	\end{array}
	$$
	Multiplying these relations, we obtain~\eqref{Qf-3x3} with $z = \boldsymbol{\O}$.
\end{proof}

\begin{proof}[Proof of Theorem~\ref{Th-Qfin->P}]
	Let ${\rm Q}^f = \{ q_V^z, z \in \widehat{X_V^f}, V \in W \}$ be an $f$--specification. For any $V \in W$ and $x \in X^V$, put
	\begin{equation}
		\label{Q->P}
		P_V (x) = \frac{q_V^y (x)}{q_I^x (y)} \left(\sum \limits_{\alpha \in X^V} \frac{q_V^y (\alpha)}{q_I^\alpha (y)} \right)^{-1},
	\end{equation}
	where $y \in X^I$, $I \in W\left(\mathbb{Z}^d \backslash V \right)$. Let us show that this formula is correct (the values of $P_V$ do not depend on the choice of $y$ and $I$), and the family $P_{Q^f} = \{P_V, V \in W \}$ is a consistent in Kolmogorov's sense system of probability distributions.
	
	First, we verify that the values of $P_V$ do not depend on the choice of $y$. From~\eqref{Qf-4x4}, it follows that for any configuration $v \in X^I$, we have
	$$
	\frac{q_V^y(x)}{q_I^x(y)} \cdot \frac{q_V^v(u)}{q_I^u(v)} = \frac{q_V^v(x)}{q_I^x(v)} \cdot \frac{q_V^y(u)}{q_I^u(y)},
	$$
	and hence,
	$$
	\frac{q_V^y(x)}{q_I^x(y)} \left( \sum \limits_{u \in X^V} \frac{q_V^y(u)}{q_I^u(y)} \right)^{-1} = \frac{q_V^v(x)}{q_I^x(v)} \left( \sum \limits_{u \in X^V} \frac{q_V^v(u)}{q_I^u(v)} \right)^{-1}.
	$$
	
	Now we show that the values of $P_V$ do not depend on the choice of $I$. Let $J \in W(\mathbb{Z}^d \backslash V)$ be another set. First, suppose that $J \cap I = {\O}$. According to~\eqref{Qf-3x3} with $z = \boldsymbol{\O}$, for any $x,\alpha \in X^V$, $y \in X^I$, $w \in X^J$, we have
	$$
	q_V^w(x) q_I^x(y) q_J^y(w) = q_I^w(y) q_V^y(x) q_J^x(w)
	$$
	and
	$$
	q_V^w(\alpha) q_I^\alpha(y) q_J^y(w) = q_I^w(y) q_V^y(\alpha) q_J^\alpha(w).
	$$
	Taking the ratio of the corresponding parts of these two equalities, we obtain
	$$
	\frac{q_V^w(x)}{q_J^x(w)} \cdot \frac{q_V^y(\alpha)}{q_I^\alpha(y)} = \frac{q_V^y(x)}{q_I^x(y)} \cdot \frac{q_V^w(\alpha)}{q_J^\alpha(w)}.
	$$
	From here it follows that
	$$
	\frac{q_V^w(x)}{q_J^x(w)} \left( \sum \limits_{\alpha \in X^V} \frac{q_V^w(\alpha)}{q_J^\alpha(w)} \right)^{-1} = \frac{q_V^y(x)}{q_I^x(y)} \left( \sum \limits_{\alpha \in X^V} \frac{q_V^y(\alpha)}{q_I^\alpha(y)} \right)^{-1}.
	$$
	Suppose now that $I \cap J = S \neq {\O}$. It is sufficient to show that for any $x, \alpha \in X^V$, $y \in X^{I \backslash S}$, $w \in X^{J \backslash S}$ and $z \in X^S$, one has
	$$
	\frac{q_V^{zw}(x)}{q_J^x (zw)} \cdot \frac{q_V^{zy}(\alpha)}{q_I^\alpha(zy)} = \frac{q_V^{zy}(x)}{q_I^x(zy)} \cdot \frac{q_V^{zw}(\alpha)}{q_J^\alpha (zw)} .
	$$
	According to~\eqref{Qf-sogl}, this relation is equivalent to
	$$
	\frac{q_V^{zw}(x)}{q_{J \backslash S}^{zx} (w)} \cdot \frac{q_V^{zy}(\alpha)}{q_{I \backslash S}^{z\alpha}(y)} = \frac{q_V^{zy}(x)}{q_{I \backslash S}^{zx}(y)} \cdot \frac{q_V^{zw}(\alpha)}{q_{J \backslash S}^{z\alpha} (w)},
	$$
	which holds true due to~\eqref{Qf-3x3}.
	
	From~\eqref{Q->P}, it obviously follows that the function $P_V$ is a probability distribution on $X^V$, $V \in W$. Further, let us verify that the system $\{P_V, V \in W\}$ is consistent in Kolmogorov's sense. Using~\eqref{Qf-sogl} and~\eqref{Q-sogl-ekv1}, for any disjoint sets $V,\Lambda \in W$ and all $x \in X^V$, we can write
	$$
	\begin{array}{l}
		\displaystyle\sum \limits_{v \in X^\Lambda} P_{V \cup \Lambda} (xv) = \sum \limits_{v \in X^\Lambda} \dfrac{q_{V \cup \Lambda}^y (xv)}{q_I^{xv} (y)} \left(\sum \limits_{\alpha \in X^V, \beta \in X^\Lambda} \dfrac{q_{V \cup \Lambda}^y (\alpha \beta)}{q_I^{\alpha \beta} (y)} \right)^{-1} = \\
		\\
		= \displaystyle\sum \limits_{v \in X^\Lambda} \dfrac{q_V^y (x) q_\Lambda^{yx} (v) q_\Lambda^x (v)}{q_I^x (y) q_\Lambda^{xy} (v)} \left(\sum \limits_{\alpha \in X^V, \beta \in X^\Lambda} \dfrac{q_V^y (\alpha) q_\Lambda^{y \alpha} (\beta) q_\Lambda^\alpha (\beta)}{q_I^\alpha (y) q_\Lambda^{\alpha y} (\beta)} \right)^{-1}= \\
		\\
		= \dfrac{q_V^y (x)}{q_I^x (y)} \displaystyle\sum \limits_{v \in X^\Lambda} q_\Lambda^x (v) \left(\sum \limits_{\alpha \in X^V} \dfrac{q_V^y (\alpha)}{q_I^\alpha (y)} \sum \limits_{\beta \in X^\Lambda} q_\Lambda^\alpha (\beta) \right)^{-1} = \\
		\\
		= \dfrac{q_V^y (x)}{q_I^x (y)} \left(\sum \limits_{\alpha \in X^V} \dfrac{q_V^y (\alpha)}{q_I^\alpha (y)} \right)^{-1} = P_V (x),
	\end{array}
	$$
	where $y \in X^I$ and $I \in W \left(\mathbb{Z}^d \backslash (V \cup \Lambda) \right)$. Thus, we showed that there exists a random field $P_{{\rm Q}^f} = \{P_V, V \in W\}$ constructed by $f$--specification ${\rm Q}^f$.
	
	Let us show that $P_{{\rm Q}^f}$ is compatible with ${\rm Q}^f$, that is, that the $f$--distri\-bu\-tion $Q^f (P_{{\rm Q}^f}) = \{Q_V^z, z \in \widehat{X_V^f}, V \in W\}$ of the random field $P_{{\rm Q}^f}$ coincides with ${\rm Q}^f$. Applying~\eqref{P->QP},~\eqref{Q->P} and~\eqref{Q-sogl-ekv1}, for any disjoint sets $V,\Lambda \in W$ and all $x \in X^V$, $z \in X^\Lambda$, we can write
	$$
	\begin{array}{l}
		Q_V^z (x) = \dfrac{P_{V \cup \Lambda}(xz)}{P_\Lambda (z)} = \\
		\\
		= \dfrac{q_{V \cup \Lambda}^y (xz)}{q_I^{xz} (y)} \left( \displaystyle\sum \limits_{\alpha \in X^V, \beta \in X^\Lambda} \dfrac{q_{V \cup \Lambda}^y (\alpha \beta)}{q_I^{\alpha \beta} (y)} \right)^{-1} \cdot \dfrac{q_I^z (y)}{q_\Lambda^y (z)} \displaystyle\sum \limits_{\beta \in X^\Lambda} \dfrac{q_\Lambda^y (\beta)}{q_I^\beta (y)} = \\
		\\
		= \dfrac{q_\Lambda^y (z) q_V^{yz}(x) q_V^z(x)}{q_I^z(y) q_V^{zy}(x)} \left( \displaystyle\sum \limits_{\alpha \in X^V, \beta \in X^\Lambda} \dfrac{q_\Lambda^y (\beta) q_V^{y \beta} (\alpha) q_V^\beta (\alpha)}{q_I^\beta (y) q_V^{\beta y} (\alpha)} \right)^{-1} \dfrac{q_I^z (y)}{q_\Lambda^y (z)} \displaystyle\sum \limits_{\beta \in X^\Lambda} \dfrac{q_\Lambda^y (\beta)}{q_I^\beta (y)} = \\
		\\
		= q_V^z(x) \left( \displaystyle\sum \limits_{\beta \in X^\Lambda} \dfrac{q_\Lambda^y (\beta)}{q_I^\beta (y) } \sum \limits_{\alpha \in X^V} q_V^\beta (\alpha) \right)^{-1} \displaystyle\sum \limits_{\beta \in X^\Lambda} \dfrac{q_\Lambda^y (\beta)}{q_I^\beta (y)} = q_V^z(x),
	\end{array}
	$$
	where $y \in X^I$ and $I \in W \left( \mathbb{Z}^d \backslash (V \cup \Lambda) \right)$.
	
	It remains to note that $P_{{\rm Q}^f}$ is a unique random field compatible with ${\rm Q}^f$. Indeed, if $\hat P$ is another random field compatible with ${\rm Q}^f$, then $Q^f_{\hat P} = {\rm Q}^f = Q^f_P$, and by Theorem~\ref{Th-QPfin->P}, $\hat P = P$. Therefore, ${\rm Q}^f$ specifies $P_{{\rm Q}^f}$.
\end{proof}

Note that formula~\eqref{Q->P} can be written in the following equivalent form
\begin{equation}
	\label{Q->P-2}
	P_V (x) = \left(\sum \limits_{\beta \in X^I} \frac{q_I^x (\beta)}{q_V^\beta (x)} \right)^{-1}, \qquad x \in X^V, V \in W,
\end{equation}
where $I \in W \left(\mathbb{Z}^d \backslash V \right)$. Indeed, for any $y \in X^I$, we can write
$$
\begin{array}{l}
	\left( \displaystyle\dfrac{q_V^y (x)}{q_I^x (y)} \left(\sum \limits_{\alpha \in X^V} \dfrac{q_V^y (\alpha)}{q_I^\alpha (y)} \right)^{-1} \right)^{-1} = \displaystyle\sum \limits_{\alpha \in X^V} \dfrac{q_I^x (y) q_V^y (\alpha)}{q_V^y (x) q_I^\alpha (y)} = \\
	\\
	= \displaystyle\sum \limits_{\alpha \in X^V} \dfrac{q_I^x (y) q_V^y (\alpha) \sum \limits_{\beta \in X^I} q_I^{\alpha}(\beta)}{q_V^y (x) q_I^\alpha (y)} = \sum \limits_{\alpha \in X^V} \sum \limits_{\beta \in X^I} \dfrac{q_I^x (y) q_V^y (\alpha) q_I^{\alpha}(\beta) q_V^{\beta}(x)}{q_V^y (x) q_I^\alpha (y) q_V^{\beta}(x)} = \\
	\\
	= \displaystyle\sum \limits_{\alpha \in X^V} \sum \limits_{\beta \in X^I} \dfrac{q_I^x (\beta) q_V^{\beta} (\alpha) q_I^{\alpha}(y) q_V^{y}(x)}{q_V^y (x) q_I^\alpha (y) q_V^{\beta}(x)} = \sum \limits_{\beta \in X^I} \dfrac{q_I^x (\beta)}{ q_V^{\beta}(x)} \sum \limits_{\alpha \in X^V} q_V^{\beta} (\alpha) = \sum \limits_{\beta \in X^I} \dfrac{q_I^x (\beta)}{ q_V^{\beta}(x)},
\end{array}
$$
where we used~\eqref{Qf-4x4} and the fact that $\sum\limits_{\beta \in X^I } {q_I^\alpha (\beta)} = \sum\limits_{\alpha \in X^V } {q_V^\beta (\alpha)} = 1$.

Note also that the results of Theorems~\ref{Th-QPfin->P} and~\ref{Th-Qfin->P} can be formulated in the following equivalent form.

\begin{theorem}
	A system ${\rm Q}^f = \{ q_V^z, z \in \widehat{X_V^f}, V \in W \}$ of strictly positive finite-dimensional distributions parameterized by finite boundary conditions is an $f$--distribution of the unique random field $P$ compatible with it if and only if the elements of ${\rm Q}^f$ satisfy the consistency conditions~\eqref{Qf-sogl}.
\end{theorem}

From the theorems above, it follows that there is a one-to-one correspondence between a random field $P$ and an $f$--specification ${\rm Q}^f$. In this regard, the random field $P$ can be identified with its system of finite-conditional distributions, and one can write $P = \{Q_V^z, z \in \widehat{X_V^f}, V \in W\}$. Therefore, there are no statements about random fields that cannot be expressed in terms of their $f$--distributions.

%

For example, in terms of ${\rm Q}^f$, estimates for mixing coefficients for random fields with weakly dependent components can be obtained. Following Dobrushin, Dala\-lyan and Nahapetian (see Theorem 2 in~\cite{DN2011}) gave an estimate for the difference between the conditional and unconditional distributions of a random field, expressed by the difference between its one-point conditional distributions with finite boundary conditions that differ at a point. Below we present another proof of this result.

\begin{prop}
	Let $P = \{Q_V^z, z \in \widehat{X_V^f}, V \in W\}$ be a random field. Then for any $V,\Lambda \in W$, $V \cap \Lambda = {\O}$, the following relation holds
	$$
	\sup \limits_{x \in X^V, z \in X^\Lambda} \left| P_V(x) - Q_V^z(x) \right| \le \sum \limits_{t \in V} \sum \limits_{s \in \Lambda} \rho_{ts},
	$$
	where
	$$
	\rho_{ts} = \sup \limits_{w \in \widehat{X^f_{\{t,s\}}}} \sup \limits_{y,v \in X^s} \sup \limits_{x \in X^t} \left| Q_t^{wy}(x) - Q_t^{wv}(x) \right|, \qquad t,s \in \mathbb{Z}^d.
	$$
\end{prop}
\begin{proof}
	For any disjoint sets $V, \Lambda \in W$ and any $x \in X^V$, $z \in X^\Lambda$, we have
	$$
	P_V(x) - Q_V^z(x) 
	= \sum \limits_{w \in X^\Lambda} P_\Lambda(w) \left( Q_V^w(x) - Q_V^z(x) \right).
	$$
	Let $\Lambda = \{s_1, s_2, ..., s_m\}$ be some enumeration of the points of $\Lambda$, $m = \vert \Lambda \vert \ge 1$. Using notations~\eqref{xu-notations}, we can write
	$$
	Q_V^w(x) - Q_V^z(x) = \sum \limits_{k=1}^m \left( Q_V^{(zw)_k w_k}(x) - Q_V^{(zw)_k z_k}(x) \right).
	$$
	Let now $V = \{t_1, t_2, ..., t_n\}$ be some enumeration of the points of $V$, $n = \vert V \vert \ge 1$. For each $k$, $1 \le k \le m$, denoting for brevity $y = (zw)_k$, $\alpha = z_k$ and $\beta = w_k$, we have
	$$
	\begin{array}{l}
		Q_V^{(zw)_k w_k}(x) - Q_V^{(zw)_k z_k}(x) = Q_V^{y \beta}(x) - Q_V^{y\alpha}(x) = \\
		\\
		= \left( Q_{t_1}^{y \beta}(x_1) - Q_{t_1}^{y \alpha}(x_1) \right)Q_{V \backslash \{t_1\}}^{y \beta x_1} (x_{V \backslash \{t_1\}}) + Q_{t_1}^{y \alpha}(x_1) \left( Q_{V \backslash \{t_1\}}^{y \alpha x_1}(x_{V \backslash \{t_1\}}) - Q_{V \backslash \{t_1\}}^{y \alpha x_1}(x_{V \backslash \{t_1\}}) \right).
	\end{array}
	$$
	Similarly, for the bracketed expression in the right-hand summand, we obtain
	$$
	\begin{array}{l}
		Q_{V \backslash \{t_1\}}^{y \alpha x_1}(x_{V \backslash \{t_1\}}) - Q_{V \backslash \{t_1\}}^{y \alpha x_1}(x_{V \backslash \{t_1\}}) = \\
		\\
		= \left( Q_{t_2}^{y x_1 \beta}(x_2) - Q_{t_2}^{y x_1 \alpha}(x_2) \right) Q_{V \backslash \{t_2\}}^{yx_1 x_2 \beta} (x_{V \backslash \{t_1,t_2\}}) 
		+ Q_{t_2}^{y \alpha x_1}(x_2) \left( Q_{V \backslash \{t_2\}}^{yx_1 x_2 \beta} (x_{V \backslash \{t_1,t_2\}}) - Q_{V \backslash \{t_2\}}^{yx_1 x_2 \alpha} (x_{V \backslash \{t_1,t_2\}}) \right).
	\end{array}
	$$
	Continuing this process the required number of times, we get
	$$
	\begin{array}{l}
		Q_V^{y \beta}(x) - Q_V^{y \alpha}(x) = \displaystyle\sum \limits_{j=1}^n Q_{t_1}^{y \alpha}(x_1) Q_{t_2}^{y \alpha x_1}(x_2) \cdot ... \cdot Q_{t_{j-1}}^{y \alpha x_1 x_2 ... x_{j-2}}(x_{j-1}) \cdot \\
		\\
		\qquad \cdot \left( Q_{t_j}^{yx_1...x_{j-1} \beta}(x_j) - Q_{t_j}^{yx_1...x_{j-1} \alpha}(x_j) \right) Q_{V \backslash \{t_1,...,t_j\}}^{y \beta x_1 ... x_j} (x_{V \backslash \{t_1,...,t_j\}}) =\\
		\\
		= \displaystyle\sum \limits_{j=1}^n Q_{\{t_1,...,t_{j-1}\}}^{y \alpha}(x_{\{t_1,...,t_{j-1}\}}) \cdot  \left( Q_{t_j}^{yx_1...x_{j-1} \beta}(x_j) - Q_{t_j}^{yx_1...x_{j-1} \alpha}(x_j) \right) Q_{\{t_{j+1},...,t_n\}}^{y \beta x_1 ... x_j} (x_{\{t_{j+1},...,t_n\}}).
	\end{array}
	$$
	
	Finally, we obtain
	$$
	\begin{array}{l}
		Q_V^w(x) - Q_V^z(x) =  \displaystyle\sum \limits_{k=1}^m \sum \limits_{j=1}^n Q_{\{t_1,...,t_{j-1}\}}^{(zw)_k z_k}(x_{\{t_1,...,t_{j-1}\}})  \cdot\\
		\\
		\quad \cdot \left( Q_{t_j}^{x_1...x_{j-1} (zw)_k w_k}(x_j) - Q_{t_j}^{x_1...x_{j-1} (zw)_k z_k}(x_j) \right) Q_{\{t_{j+1},...,t_n\}}^{x_1 ... x_j (zw)_k w_k} (x_{\{t_{j+1},...,t_n\}}).
	\end{array}
	$$
	From here it follows that
	$$
	\begin{array}{l}
		\sup \limits_{x \in X^V, z \in X^\Lambda} \left| P_V(x) - Q_V^z(x) \right| \le  \sup \limits_{x \in X^V, z \in X^\Lambda} \displaystyle\sum \limits_{w \in X^\Lambda} P_\Lambda(w) \sum \limits_{k=1}^m \sum \limits_{j=1}^n Q_{\{t_1,...,t_{j-1}\}}^{(zw)_k z_k}(x_{\{t_1,...,t_{j-1}\}}) \cdot \\
		\\
		\quad \cdot  \left| Q_{t_j}^{x_1...x_{j-1} (zw)_k w_k}(x_j) - Q_{t_j}^{x_1...x_{j-1} (zw)_k z_k}(x_j) \right| Q_{\{t_{j+1},...,t_n\}}^{x_1 ... x_j (zw)_k w_k} (x_{\{t_{j+1},...,t_n\}}) \le \\
		\\
		\le \displaystyle\sum \limits_{k=1}^m \sum \limits_{j=1}^n \rho_{t_j s_k} = \sum \limits_{t \in V} \sum \limits_{s \in \Lambda} \rho_{ts}.
	\end{array}
	$$
\end{proof}

Further, it seems more natural to give the definition of a Markov random field in terms of the elements of its $f$--distribution (see, for example,~\cite{Griff}). Namely, a random field $P$ will be called a Markov random field (with respect to a neighborhood system $\partial$ on $\mathbb{Z}^d$) if the elements of its $f$--distribution $Q_P^f$ satisfy the following conditions: for all $V \in W$ and $z \in \widehat{X_V^f}$ such that $\partial V \subset s(z)$, it holds
\begin{equation}
	\label{MarkovV}
	Q_V^z (x) = Q_V^{z_{\partial V}} (x), \qquad x \in X^V.
\end{equation}
Note that Dobrushin~\cite{D68} defined a Markov random field somewhat differently, imposing restrictions on its conditional probabilities with infinite boundary conditions. In Section~\ref{Section-QD}, we will show the equivalence of these definitions.

\begin{remark}
	\label{remark-Q-sogl-Kolm}
	The collection of the elements of ${\rm Q}^f$ with the same boundary condition is consistent in Kolmogorov's sense. Namely, for fixed $\Lambda \in W$ and $z \in X^\Lambda$, it holds
	\begin{equation}
		\label{Q-sogl-Kolm}
		\sum \limits_{y \in X^I} q_{V \cup I}^z (xy) = q_V^z(x), \qquad x \in X^V,
	\end{equation}
	where $V,I\in W(\mathbb{Z}^d \backslash \Lambda)$, $V \cap I = {\O}$. This means that the system ${\rm Q}^{\Lambda, z} = \{q_V^z, V \in W(\mathbb{Z}^d \backslash \Lambda)\}$ of probability distributions defines a unique random field  $P^{\Lambda, z}$ on $(X^{\mathbb{Z}^d \backslash \Lambda}, \mathscr{B}^{\mathbb{Z}^d \backslash \Lambda})$.
\end{remark}

\begin{remark}
	\label{R2}
	It follows directly from~\eqref{Qf-sogl} that for the elements of ${\rm Q}^f$, it holds
	\begin{equation}
		\label{Q-sogl-ekv}
		q_{V \cup I}^z (xy) q_I^{zx} (v) = q_{V \cup I}^z (xv) q_I^{zx} (y),
	\end{equation}
	where $x \in X^V$, $y,v \in X^I$, $z \in \widehat{X_{V \cup I}^f}$, $V,I \in W$, $V \cap I = {\O}$.
	
	On the other hand, if the elements of some system of strictly positive finite-dimensional distributions parameterized by finite boundary conditions satisfy conditions~\eqref{Q-sogl-Kolm} and~\eqref{Q-sogl-ekv}, then they satisfy conditions~\eqref{Qf-sogl} as well. To verify this, it is enough to take a sum of both sides of~\eqref{Q-sogl-ekv} over all $v \in X^I$.
	
	Relations~\eqref{Q-sogl-ekv}, in their turn, hold if and only if for any $V \in W$, $s \in \mathbb{Z}^d \backslash V$ and $x,u \in X^V$, $y \in X^s$, $z \in \widehat{X_{V \cup \{s\}}^f}$, the following equality takes place
	\begin{equation}
		\label{Q-sogl-ekv-1}
		q_{V \cup \{s\}}^z (xy) q_V^{zy} (u) = q_{V \cup \{s\}}^z (uy) q_V^{zy} (x).
	\end{equation}
\end{remark}

\begin{remark}
	\label{Remark-z-zy}
	The elements of ${\rm Q}^f$ for any $V \in W$ and $z \in \widehat{X_V^f}$, satisfy the following relations
	\begin{equation}
		\label{qz-qzy}
		q_V^z (x) = \frac{q_V^{zy}(x)}{q_I^{zx}(y)} \left(\sum \limits_{\alpha \in X^V} \frac{q_V^{zy}(\alpha)}{q_I^{z \alpha}(y)} \right)^{-1}, \qquad x \in X^V,
	\end{equation}
	where $y \in X^I$, $I \in W\left(\mathbb{Z}^d \backslash (V \cup s(z)) \right)$. Note that for $z = \boldsymbol{\O}$, these relations lead to~\eqref{Q->P}.
\end{remark}

Since $f$--distribution $Q_P^f$ of any random field $P$ forms an $f$--specification, all the above remarks stay true for $Q_P^f$.

\bigskip

The connection between unconditional and finite-conditional distributions of a random field $P$ can be also expressed in the following form:
$$
P_V(x) = \sum \limits_{z \in X^\Lambda} Q_V^z(x) P_\Lambda(z), \qquad x \in X^V,
$$
where $V, \Lambda \in W$, $V \cap \Lambda = {\O}$. These relations can be considered as a finite-dimensional analogue of the well-known DLR--equations (named after Dobrushin, Lanford and Ruelle) in statistical physics. It is easy to see that the solution to the direct problem given by formula~\eqref{QP->P}, when substituted into this equation, leads to an identity. In terms of DLR--equations, the solution to the inverse problem for $Q^f$ can be stated as follows.

\begin{theorem}
	Let ${\rm Q}^f$ be an $f$--specification. Then there exists a unique random field $P$ satisfying the finite-volume DLR-equations
	\begin{equation}
		\label{DLR-fin}
		P_V(x) = \sum \limits_{z \in X^\Lambda} q_V^z(x) P_\Lambda(z), \qquad x \in X^V, V, \Lambda \in W, V \cap \Lambda = {\O}.
	\end{equation}
	In this case, $Q^f_P = {\rm Q}^f$.
\end{theorem}
\begin{proof}
	Let us show that the functions defined by~\eqref{Q->P} form a solution to  equations~\eqref{DLR-fin}. Taking into account~\eqref{Q->P} and~\eqref{Q->P-2}, for any pairwise disjoint sets $V,I,\Lambda \in W$, we obtain
	$$
	\begin{array}{l}
		\displaystyle\sum \limits_{z \in X^\Lambda} q_V^z (x) P_\Lambda (z) = \sum \limits_{z \in X^\Lambda} q_V^z (x) \dfrac{q_\Lambda^x (z)}{q_V^z (x)} \cdot \left(\sum \limits_{\alpha \in X^\Lambda} \dfrac{q_\Lambda^x (\alpha)}{q_V^\alpha (x)}\right)^{-1} 
		= \left(\displaystyle\sum \limits_{\alpha \in X^\Lambda} \dfrac{q_\Lambda^x (\alpha)}{q_V^\alpha (x)}\right)^{-1} = P_V(x).
	\end{array}
	$$
	
	Therefore, there exists a random field $P$ finite-dimensional distributions of which are defined by formula~\eqref{Q->P}. As it was shown in the proof of Theorem~\ref{Th-Qfin->P}, $Q^f_P = {\rm Q}^f$. Finally, since any random field is uniquely determined by its finite-conditional distribution, the random field $P$ is unique.
\end{proof}

\subsection{System of one-point distributions with finite boundary conditions}
\label{Section-Qf-1}

For a random field $P$, the set $Q_1^f(P) = \{Q_t^z, z \in \widehat{X^f_t}, t \in \mathbb{Z}^d\}$ of one-point conditional probabilities with finite boundary conditions will be called the \emph{one-point finite-conditional distribution of the random field $P$}, or, in short, its \emph{$1f$--distribution of $P$}. The system $Q_1^f(P)$ was introduced in~\cite{DN2011}.

Let us consider the direct problem for a random field $P$ in terms of the system $Q_1^f(P)$.

\begin{theorem}
	\label{Th-QPfin1->P}
	Any random field is restored by its $1f$--distribution.
\end{theorem}
\begin{proof}
	It is sufficient to note that the finite-dimensional distributions $\{P_V, \linebreak V \in W\}$ of the random field $P$ and its one-point finite-conditional probabilities $Q_1^f(P) = \{Q_t^z, z \in \widehat{X^f_t}, t \in \mathbb{Z}^d \}$ are connected in the following way
	$$
	P_V (x) = \left(\sum \limits_{y \in X^s} \frac{Q_s^{x_{t_1}} (y)}{Q_{t_1}^y (x_{t_1})} \right)^{-1} 
	Q_{t_2}^{x_{t_1}} (x_{t_2}) \cdot Q_{t_3}^{x_{t_1} x_{t_2}} (x_{t_3}) \cdot ... \cdot Q_{t_{\vert V \vert}}^{x_{t_1} ... x_{t{{\vert V \vert}-1}}} (x_{t_{\vert V \vert}}),
	$$
	where $x \in X^V$, $y \in X^s$, $s \in \mathbb{Z}^d \backslash V$, and $V = \{t_1, t_2, ..., t_{\vert V \vert}\}$ is some enumeration of the points of~${V \in W}$.
\end{proof}

The inverse problem for a system of one-point distributions parameterized by finite boundary conditions was first considered in~\cite{DN2011}. Below we slightly improve the result of~\cite{DN2011}.

As the main characterizing property of the elements of the system $Q_1^f(P)$, we consider the following easily verifiable property: for all $t,s \in \mathbb{Z}^d$ and $x \in X^t$, $y \in X^s$, $z \in \widehat{X^f_{\{t,s\}}}$, it holds
\begin{equation}
	\label{cond1-prop-main-2x2}
	Q_t^z (x) Q_s^{z x} (y) = Q_s^z (y) Q_t^{z y} (x).
\end{equation}

A system ${\rm Q}_1^f = \{q_t^z, z \in \widehat{X^f_t}, t \in \mathbb{Z}^d\}$ of strictly positive one-point probability distributions $q_t^z$ on $X^t$ parameterized by finite boundary conditions $z \in \widehat{X^f_t}$, $t \in \mathbb{Z}^d$, will be called a \emph{1--specification with finite boundary conditions} (or \emph{$1f$--specification}) if its elements satisfy the following consistency conditions: for all $t,s \in \mathbb{Z}^d$ and $x \in X^t$, $y \in X^s$, $z \in \widehat{X^f_{\{t,s\}}}$, it holds
\begin{equation}
	\label{Q1f-sogl-2x2}
	q_t^z (x) q_s^{zx} (y) = q_s^z (y) q_t^{zy} (x).
\end{equation}

\begin{theorem}
	\label{T-Q1f-P}
	Any $1f$--specification specifies compatible with it random field.
\end{theorem}

In the proof of this result, we need the following properties of the elements of ${\rm Q}_1^f$.

\begin{lemma}
	\label{Lemma-Q1f}
	The elements of $1f$--specification ${\rm Q}_1^f$ satisfy the following relation
	\begin{equation}
		\label{Q1f-4x4}
		q_t^{zy} (x) q_s^{zx} (v) q_t^{zv} (u) q_s^{zu} (y) = q_t^{zy} (u) q_s^{zu} (v) q_t^{zv} (x) q_s^{zx} (y)
	\end{equation}
	for all $t,s \in \mathbb{Z}^d$ and $x,u \in X^t$, $y,v \in X^s$, $z \in \widehat{X^{f0}_{\{t,s\}}}$. Further, for any points $t,s,r \in \mathbb{Z}^d$ and any configurations $x \in X^t$, $y \in X^s$, $w \in X^r$, one has
	\begin{equation}
		\label{Q1f-3x3}
		q_t^w(x) q_s^x(y) q_r^y(w) = q_s^w(y) q_t^y(x) q_r^x(w) .
	\end{equation}
\end{lemma}
\begin{proof}
	For any $t,s \in \mathbb{Z}^d$ and $x,u \in X^t$, $y,v \in X^s$, $z \in \widehat{X^f_{\{t,s\}}}$, by~\eqref{Q1f-sogl-2x2}, we can write
	\begin{multline*}
		q_s^z(y) q_t^{zy}(x) \cdot q_t^z(x) q_s^{zx}(v) \cdot q_s^z(v) q_t^{zv}(u) \cdot q_t^z(u) q_s^{zu}(y) = \\
		= q_t^z(x) q_s^{zx}(y) \cdot q_s^z(v) q_t^{zv}(x) \cdot q_t^z(u) q_s^{zu}(v) \cdot q_s^z(y) q_t^{zy}(u),
	\end{multline*}
	whence, after the necessary reductions, follows~\eqref{Q1f-4x4}. Now let us show that~\eqref{Q1f-4x4} stays true for $z = \boldsymbol{{\O}}$ as well. For different points $t,s,r \in \mathbb{Z}^d$ and any configurations $x,u \in X^t$, $y,v \in X^s$, $z \in X^r$, according to~\eqref{Q1f-sogl-2x2}, we have
	\begin{multline*}
		q_t^y(x) q_r^{xy}(z) \cdot q_s^x(v) q_r^{xv}(z) \cdot q_t^v(u) q_r^{uv}(z) \cdot q_s^u(y) q_r^{uy}(z) =\\
		= q_r^y(z) q_t^{zy}(x) \cdot q_r^x(z) q_s^{zx}(v) \cdot q_r^v(z) q_t^{zv}(u) \cdot q_r^u(z) q_s^{zu}(y)
	\end{multline*}
	and
	\begin{multline*}
		q_s^x(y) q_r^{xy}(z) \cdot q_t^y(u) q_r^{uy}(z) \cdot q_s^u(v) q_r^{uv}(z) \cdot q_t^v(x) q_r^{xv}(z) =\\
		=q_r^x(z) q_s^{zx}(y) \cdot q_r^y(z) q_t^{zy}(u) \cdot q_r^u(z) q_s^{zu}(v) \cdot q_r^v(z) q_t^{zv}(x).
	\end{multline*}
	Dividing the first of these equalities by the second one, we obtain
	$$
	\frac{q_t^y(x) q_s^x(v) q_t^v(u) q_s^u(y)}{q_s^x(y) q_t^y(u) q_s^u(v) q_t^v(x)} =
	\frac{q_t^{zy}(x) q_s^{zx}(v) q_t^{zv}(u) q_s^{zu}(y)}{q_s^{zx}(y) q_t^{zy}(u) q_s^{zu}(v) q_t^{zv}(x)}.
	$$
	It remains to note that the right-hand side of the relation above is equal to one.
	
	The validity of~\eqref{Q1f-3x3} can be shown using the same reasonings that we used to verify~\eqref{Qf-3x3}.
\end{proof}

\begin{proof}[Proof of Theorem~\ref{T-Q1f-P}]
	Let ${\rm Q}_1^f = \{q_t^z, z \in \widehat{X_t^f}, t \in \mathbb{Z}^d\}$ be a $1f$--specification. For all $V \in W$ and $x \in X^V$, put
	\begin{equation}
		\label{Q1->P}
		P_V (x) = P_{t_1} (x_{t_1}) q_{t_2}^{x_{t_1}} (x_{t_2}) \cdot q_{t_3}^{x_{t_1} x_{t_2}} (x_{t_3}) \cdot ... \cdot q_{t_{\vert V \vert}}^{x_{t_1} ... x_{t{{\vert V \vert}-1}}} (x_{t_{\vert V \vert}}),
	\end{equation}
	where
	\begin{equation}
		\label{Q1->P1}
		P_t(u) = \frac{q_t^y(u)}{q_s^u(y)} \left( \sum \limits_{\alpha \in X^t} \frac{q_t^y(\alpha)}{q_s^\alpha(y)} \right)^{-1}, \qquad u \in X^t,
	\end{equation}
	$y \in X^s$, $s \in \mathbb{Z}^d \backslash V$, and $V = \{t_1, t_2, ..., t_{\vert V \vert}\}$ is some enumeration of the points of~$V$. First, let us verify the correctness of these formulas.
	
	Using the same reasoning as in the proof of Theorem~\ref{Th-Qfin->P} and relations~\eqref{Q1f-4x4} and~\eqref{Q1f-3x3}, one can verify that the right-hand side of~\eqref{Q1->P1} does not depend on the choice of $y \in X^s$ and $s \in \mathbb{Z}^d \backslash V$. Further, note (see the derivation of formula~\eqref{Q->P-2}) that
	\begin{equation}
		\label{Q1->P1-2}
		P_t(u) = \left( \sum \limits_{\beta \in X^s} \frac{q_s^u(\beta)}{q_t^\beta(u)} \right)^{-1}, \qquad u \in X^t.
	\end{equation}
	
	Let us show that the right-hand side of~\eqref{Q1->P} does not depend on the enumeration of the elements of $V$. It is sufficient to verify this statement for enumerations $t_1,...,t_{k-1}, t_{k}, ...,t_{\vert V \vert}$ and $t_1,...,t_{k}, t_{k-1}, ...,t_{\vert V \vert}$ differing in the position of two successive points $t_{k-1}$ and $t_k$, $2 \le k \le {\vert V \vert}$. Thus, we need to check that the following equalities hold true:
	$$
	q_{t_{k - 1}}^{ x_{t_1} ...x_{t_{k - 2}} } (x_{t_{k - 1}}) q_{t_k}^{ x_{t_1} ...x_{t_{k - 2}} x_{t_{k - 1}} } (x_{t_k}) = q_{t_k}^{ x_{t_1} ... x_{t_{k - 2}} } (x_{t_k}) q_{t_{k - 1}}^{ x_{t_1} ... x_{t_{k - 2}} x_{t_k} } (x_{t_{k - 1}})
	$$
	and
	$$
	P_{t_1} (x_{t_1}) q_{t_2}^{ x_{t_1} } (x_{t_2}) = P_{t_2} (x_{t_2}) q_{t_1}^{ x_{t_2} } (x_{t_1}).
	$$
	The first relation directly follows from the consistency conditions~\eqref{Q1f-sogl-2x2}. The second one becomes obvious if we use formula~\eqref{Q1->P1-2} with $s = t_2$ and $y=x_{t_2}$ to express $P_{t_1}$ and formula~\eqref{Q1->P1} with $s = t_1$ and $y = x_{t_1}$ for $P_{t_2}$.
	
	It is not difficult to see that for each $V \in W$, the function $P_V$ defined by~\eqref{Q1->P} is a probability distribution on $X^V$, and the system $\{ P_V, V \in W \}$ is consistent in Kolmogorov's sense. Thus, there exists a random field $P_{Q_1^f}$ such that $\left(P_{{\rm Q}_1^f}\right)_V = P_V$, $V \in W$. For this random field, for any $t \in \mathbb{Z}^d$ and $z \in \widehat{X^f_t}$, we have
	$$
	Q_t^z(x) = \frac{P_{\{t\} \cup s(z)}(xz)}{P_{s(z)} (z)} = \frac{P_{s(z)} (z) q_t^z(x)}{P_{s(z)} (z)} = q_t^z(x), \qquad x \in X^t,
	$$
	and hence, $Q_1^f(P_{{\rm Q}_1^f}) = {\rm Q}_1^f$. According to Theorem~\ref{Th-QPfin1->P}, $P_{{\rm Q}_1^f}$ is the unique random field compatible with ${\rm Q}_1^f$.
\end{proof}

Note that in Theorem 1 in~\cite{DN2011}, the conditions~\eqref{Q1f-sogl-2x2} together with the relations~\eqref{Q1f-4x4} with $z = \boldsymbol{{\O}}$ were considered as consistency conditions for the system ${\rm Q}_1^f$. However, as we have seen in Lemma~\ref{Lemma-Q1f}, the relations~\eqref{Q1f-4x4} follow from~\eqref{Q1f-sogl-2x2}, and thus, for the existence of the unique random field $P_{{\rm Q}_1^f}$ it is necessary and sufficient to require the fulfilment only of the consistency conditions~\eqref{Q1f-sogl-2x2}.

From the theorems above, it follows that there is a one-to-one correspondence between a random field $P$ and an $1f$--specification ${\rm Q}_1^f$. In this regard, the random field $P$ can be identified with its system of one-point finite-conditional distributions, and one can write $P = \{Q_t^z, z \in \widehat{X^f_t}, t \in \mathbb{Z}^d\}$.

Now let us consider the relation between the systems ${\rm Q}_1^f$ and ${\rm Q}^f$. The following statement holds true.

\begin{theorem}
	\label{Th-Q1f->Qf}
	A set ${\rm Q}_1^f = \{q_t^z, z \in \widehat{X^f_t}, t \in \mathbb{Z}^d\}$ of strictly positive one-point probability distributions parameterized by finite boundary conditions is a one-point subsystem of an $f$--specification ${\rm Q}^f = \{q_V^z, z \in \widehat{X^f_V}, V \in W\}$ if and only if ${\rm Q}_1^f$ is a $1f$--specification. The specification ${\rm Q}^f$ is uniquely determined by~${\rm Q}_1^f$.
\end{theorem}
\begin{proof}
	The necessity follows from the fact that the consistency conditions~\eqref{Q1f-sogl-2x2} of the elements of $1f$--specification ${\rm Q}_1^f$ coincide with the property~\eqref{Q-sogl-ekv1} of the elements of $f$--specification ${\rm Q}^f$ for $V = \{t\}$ and $I = \{s\}$. Let us prove the sufficiency.
	
	Let ${\rm Q}_1^f = \{q_t^z, z \in \widehat{X^f_t}, t \in \mathbb{Z}^d\}$ be a $1f$--specification. For any $V \in W$, put
	\begin{equation}
		\label{Q1->Q-1}
		q_V^z(x) = q_{t_1}^z(x_{t_1}) q_{t_2}^{zx_{t_1}}(x_{t_2}) \cdot ... \cdot q_{t_{\vert V \vert}}^{z x_{t_1} x_{t_2} ... x_{t_{\vert V \vert -1}}}(x_{t_{\vert V \vert}}),
	\end{equation}
	where $V = \{t_1, t_2, ..., t_{\vert V \vert}\}$ is some enumeration of the points of the set $V$. By~\eqref{Q1f-sogl-2x2}, the values of $q_V^z$ do not depend on the way of enumeration of the points of $V$. Moreover, as it is easy to see, for each $V \in W$ and $z \in \widehat{X_V^f}$, the function $q_V^z$ is a positive probability distribution on $X^V$.
	
	Let us show that the system ${\rm Q}^f = \{q_V^z, z \in \widehat{X_V^f}, V \in W\}$ forms an $f$--specification. For any disjoint sets $V = \{t_1, t_2, ..., t_{\vert V \vert}\}, I = \{s_1, s_2,...,s_{\vert I \vert}\} \in W$ and configurations $x \in X^V$, $y \in X^I$, $z \in \widehat{X^f_{V \cup I}}$, one has
	$$
	\begin{array}{l}
		q_{V \cup I}^z(xy) = q_{t_1}^z(x_{t_1}) q_{t_2}^{zx_{t_1}}(x_{t_2}) \cdot ... \cdot q_{t_{\vert V \vert}}^{z x_{t_1} x_{t_2} ... x_{t_{\vert V \vert -1}}}(x_{t_{\vert V \vert}}) \cdot \\
		\\
		\cdot q_{s_1}^{zx} (y_{s_1}) q_{s_2}^{zx y_{s_1}}(y_{s_2}) \cdot ... \cdot q_{s_{\vert I \vert}}^{z x y_{s_1} y_{s_2} ... y_{s_{\vert I \vert -1}}}(y_{s_{\vert I \vert}}) = q_V^z(x) q_I^{zx}(y).
	\end{array}
	$$
	Hence, the elements of ${\rm Q}^f$ satisfy the consistency conditions~\eqref{Qf-sogl}.
	
	The uniqueness of ${\rm Q}^f$ follows from the construction of its elements.
\end{proof}

The elements of an $f$--specification ${\rm Q}^f$ can be also constructed by the elements of $1f$--specification ${\rm Q}_1^f$ in the following way first introduced in~\cite{DN2004} for Dobrushin-type specifications. Note that this approach can be also used under less restrictive conditions than strict positivity of the elements of specifications.

\begin{prop}
	\label{Prop-Q1f->Qf}
	Let ${\rm Q}_1^f = \{q_t^z, z \in \widehat{X^f_t}, t \in \mathbb{Z}^d\}$ be a $1f$--specification.
	For any $V \in W$ and $x \in X^V$, $z \in \widehat{X_V^f}$, put
	\begin{equation}
		\label{Q1->Q-2}
		q_V^z(x) = \prod \limits_{j=1}^n \frac{q_{t_j}^{z (xu)_j} (x_j)}{q_{t_j}^{z (xu)_j} (u_j)} \cdot \left( \sum \limits_{\alpha \in X^V} \prod \limits_{j=1}^n \frac{q_{t_j}^{z (\alpha u)_j} (\alpha_j)}{q_{t_j}^{z (\alpha u)_j} (u_j)} \right)^{-1},
	\end{equation}
	where $u \in X^V$, $V = \{t_1, t_2, ..., t_n\}$ is some enumeration of the points of the set $V$, $\vert V \vert = n$, and we used notations~\eqref{xu-notations}. Then the system ${\rm Q}^f = \{q_V^z, z \in \widehat{X_V^f}, V \in W\}$ forms an $f$--specification.
\end{prop}

Formula~\eqref{Q1->Q-2} can be equivalently written in the following recurrent form
$$
q_{\{s\} \cup V}^z(yx) = \frac{q_s^{zu} (y) q_V^{zy}(u)}{q_s^{zu}(v) q_V^{zy}(x)} \left( \sum \limits_{\alpha \in X^V, \beta \in X^s} \frac{q_s^{zu} (\beta) q_V^{z\beta}(\alpha)}{q_s^{zu}(v) q_V^{z\beta}(u)} \right)^{-1},
$$
$x,u \in X^V$, $y,v \in X^s$, $z \in \widehat{X^f_{V \cup \{s\}}}$, $s \in \mathbb{Z}^d \backslash V$, $V \in W$.

The correctness of formula~\eqref{Q1->Q-2} can be checked in the same way as it was done in the paper~\cite{DN2004} (see also~\cite{KhN22}). The equivalence of formulas~\eqref{Q1->Q-1} and~\eqref{Q1->Q-2} derives from the following chain of equalities
$$
\begin{array}{l}
	q_{\{s\} \cup V}^z(yx) = \dfrac{q_s^{zu} (y) q_V^{zy}(x)}{q_s^{zu}(v) q_V^{zy}(u)} \left( \displaystyle\sum \limits_{\alpha \in X^V, \beta \in X^s} \dfrac{q_s^{zu} (\beta) q_V^{z\beta}(\alpha)}{q_s^{zu}(v) q_V^{z\beta}(u)} \right)^{-1} = \\
	\\
	= \dfrac{q_s^{zu} (y)}{q_V^{zy}(u)} \left( \displaystyle\sum \limits_{\beta \in X^s} \dfrac{q_s^{zu} (\beta)}{q_V^{z\beta}(u)} \right)^{-1} q_V^{zy}(x) = q_s^z(y) q_V^{zy}(x).
\end{array}
$$
were we used~\eqref{qz-qzy}.


Constructed by ${\rm Q}_1^f$ $f$--specification ${\rm Q}^f$ inherits properties of ${\rm Q}_1^f$, for example, Morkov property. We say that an $f$--specification ${\rm Q}^f$ is Markov if its elements satisfy the Markov property~\eqref{MarkovV}. The following statement holds true.


\begin{prop}
	For an $f$--specification ${\rm Q}^f$ to be Markov it is necessary and sufficient that its one-point subsystem ${\rm Q}_1^f$ satisfy the Markov property.
\end{prop}
\begin{proof}
	The necessity is obvious. Let us prove the sufficiency. Let the elements of ${\rm Q}_1^f$ satisfy the Markov property: for all $t \in \mathbb{Z}^d$ and $z \in \widehat{X^f_t}$ such that $\partial t \subset s(z)$, it holds
	$$
	q_t^z (x) = q_t^{z_{\partial t}} (x), \qquad x \in X^t,
	$$
	where $\partial = \{ \partial t, t \in \mathbb{Z}^d \}$ is a neighborhood system in $\mathbb{Z}^d$. For $V \in W$, let $\Lambda \in W(\mathbb{Z}^d \backslash V)$ be such that $\partial V \subset \Lambda$. Since for any $t \in V$, one has
	$$
	\partial t \subset \partial V \cup (V \backslash \{t\}) \subset \Lambda \cup (V \backslash \{t\}),
	$$
	for all $z \in X^\Lambda$ and $x,u \in X^V$, we can write
	$$
	\displaystyle\prod \limits_{j=1}^n \frac{q_{t_j}^{z_{\partial V} (xu)_j} (x_j)}{q_{t_j}^{z_{\partial V} (xu)_j} (u_j)} = \prod \limits_{j=1}^n \frac{q_{t_j}^{(z_{\partial V} (xu)_j)_{\partial t_j}} (x_j)}{q_{t_j}^{(z_{\partial V} (xu)_j)_{\partial t_j}} (u_j)} = \prod \limits_{j=1}^n \frac{q_{t_j}^{(z_\Lambda (xu)_j)_{\partial t_j}} (x_j)}{q_{t_j}^{(z_\Lambda (xu)_j)_{\partial t_j}} (u_j)} 
		=\displaystyle\prod \limits_{j=1}^n \frac{q_{t_j}^{z (xu)_j} (x_j)}{q_{t_j}^{z (xu)_j} (u_j)},
	$$
	where $V = \{t_1, t_2, ..., t_n\}$ is some enumeration of the points of $V$, $n = \vert V \vert$. Hence, using~\eqref{Q1->Q-2} and the Markov property of the elements of ${\rm Q}_1^f$, we obtain
	$$
	\begin{array}{l}
		q_V^z (x) = \displaystyle\prod \limits_{j=1}^n \dfrac{q_{t_j}^{z (xu)_j} (x_j)}{q_{t_j}^{z (xu)_j} (u_j)} \cdot \left( \sum \limits_{\alpha \in X^V} \prod \limits_{j=1}^n \dfrac{q_{t_j}^{z (\alpha u)_j} (\alpha_j)}{q_{t_j}^{z (\alpha u)_j} (u_j)} \right)^{-1} = \\
		\\
		= \displaystyle\prod \limits_{j=1}^n \dfrac{q_{t_j}^{z_{\partial V} (xu)_j} (x_j)}{q_{t_j}^{z_{\partial V} (xu)_j} (u_j)} \cdot \left( \sum \limits_{\alpha \in X^V} \prod \limits_{j=1}^n \dfrac{q_{t_j}^{z_{\partial V} (\alpha u)_j} (\alpha_j)}{q_{t_j}^{z_{\partial V} (\alpha u)_j} (u_j)} \right)^{-1} = q_V^{z_{\partial V}} (x).
	\end{array}
	$$
\end{proof}

Moreover, the following statement is true.

\begin{prop}
	Let ${\rm Q}_1^f$ be a $1f$--specification and let ${\rm Q}^f$ be the constructed by it $f$--specification. Then $P_{{\rm Q}_1^f} = P_{{\rm Q}^f}$.
\end{prop}
\begin{proof}
	It is sufficient to note that $f$--specification ${\rm Q}^f$ is uniquely defined by $1f$--specification ${\rm Q}_1^f$, while each of the systems ${\rm Q}_1^f$ and ${\rm Q}^f$ specifies a compatible with it random field.
\end{proof}

\subsection{System of Palm-type distributions}
\label{Section-Q-Pi}

Let $P = \{P_V, V \in W\}$ be a random field and let $Q^{\Pi}_P = \{ Q_V^z, z \in X^t, t \in \mathbb{Z}^d \backslash V, V \in W \}$ be a family of its conditional probabilities under the condition at a point defined by~\eqref{P->QP} for $S = \{t\}$, $t \in \mathbb{Z}^d$. The system $Q^{\Pi}_P$ will be called the \emph{Palm distribution of the random field $P$}.

\begin{theorem}
	\label{Th-QPPalm->P}
	Any random field is restored by its Palm distribution.
\end{theorem}
\begin{proof}
	It is sufficient to note that for any $V \in W$ and $x \in X^V$, the following relations hold
	$$
	P_V(x) = \left(\sum \limits_{y \in X^s} \frac{Q_s^{x_t} (y)}{Q_t^y (x_t)} \right)^{-1} Q_{V \backslash \{t\}}^{x_t} \left(x_{V \backslash \{t\}}\right),
	$$
	where $t \in V$, $s \in \mathbb{Z}^d \backslash V$.
\end{proof}

It is clear that under a fixed boundary condition, conditional probabilities are consistent in Kolmogorov's sense. Moreover,
$$
Q_t^y (x) Q_{s \cup V}^x (yu) = Q_s^x (y) Q_{t \cup V}^y (xu),
$$
where $t,s \in \mathbb{Z}^d$, $V \in W(\mathbb{Z}^d \backslash \{t,s\})$ and $x \in X^t$, $y \in X^s$, $u \in X^V$. It is these two relations between the Palm-type conditional probabilities that we propose to consider as consistency conditions for the elements of a system of distributions parameterized by a boundary condition at a point.

A set ${\rm Q}^{\Pi} = \left\{ q_V^z, z \in X^t, t \in \mathbb{Z}^d \backslash V, V \in W \right\}$ of strictly positive probability distributions $q_V^z$ on $X^V$ parameterized by boundary condition $z \in X^t$ at a point $t$, $t \in \mathbb{Z}^d$ will be called a \emph{Palm specification} if its elements satisfy the following consistency conditions:\\
1. for all $t,s \in \mathbb{Z}^d$, $V \in W(\mathbb{Z}^d \backslash \{t,s\})$ and $x \in X^t$, $y \in X^s$, $u \in X^V$, it holds
\begin{equation}
	\label{QPi-sogl}
	q_t^y (x) q_{s \cup V}^x (yu) = q_s^x (y) q_{t \cup V}^y (xu);
\end{equation}
2. for all disjoint sets $I, V \in W$, $t \in \mathbb{Z}^d \backslash (V \cup I)$ and $x \in X^V$, $z \in X^t$, it holds
\begin{equation}
	\label{QPi-sogl-I-V}
	\sum \limits_{y \in X^I} q_{V \cup I}^z (xy) = q_V^z (x).
\end{equation}

The solution to the inverse problem for the system ${\rm Q}^{\Pi}$ is given in the next theorem.

\begin{theorem}
	\label{T-QPi->P}
	Any Palm specification specifies compatible with it random field.
\end{theorem}
\begin{proof}
	Let ${\rm Q}^{\Pi} = \left\{ q_V^z, z \in X^t, t \in \mathbb{Z}^d, V \in W(\mathbb{Z}^d \backslash \{t\}) \right\}$ be a Palm specification. For any $t \in \mathbb{Z}^d$, define  $P_t$ by formula~\eqref{Q1->P1}, and for any $V \in W$, $\vert V \vert >1$, put
	\begin{equation}
		\label{QPi->P}
		P_V (x) = P_t(x_t) q_{V\backslash \{t\}}^{x_t}(x_{V\backslash \{t\}}), \qquad x \in X^V,
	\end{equation}
	where $t \in V$.
	
	First, let us verify the correctness of these formulas. By~\eqref{QPi-sogl}, for all $t,s \in \mathbb{Z}^d$, $V \in W(\mathbb{Z}^d \backslash \{t,s\})$ and $x,u \in X^t$, $y,v \in X^s$, $z \in X^V$, we have
	$$
	\begin{array}{l}
		q_t^y(x) q_{s \cup V}^x (yz) \cdot q_s^x(v) q_{t \cup V}^v(xz) \cdot q_t^v(u) q_{s \cup V}^u (vz) \cdot q_s^u(y) q_{t \cup V}^y(uz) = \\
		\\
		= q_s^x(y) q_{t \cup V}^y(xz) \cdot q_t^v(x) q_{s \cup V}^x (vz) \cdot q_s^u(v) q_{t \cup V}^v(uz) \cdot q_t^y(u) q_{s \cup V}^u (yz),
	\end{array}
	$$
	and thus,
	$$
	q_t^y (x) q_s^x (v) q_t^v (u) q_s^u (y) = q_t^y (u) q_s^u (v) q_t^v (x) q_s^x (y).
	$$
	Moreover, for any $t,s,r \in \mathbb{Z}^d$ and $x \in X^t$, $y \in X^s$, $v \in X^r$, we can write
	$$
		q_t^y (x) q_{\{s,r\}}^x (yv) \cdot q_r^x (v) q_{\{t,s\}}^v (xy) \cdot q_s^v (y) q_{\{t,r\}}^y (xv) 
		= q_s^x (y) q_{\{t,r\}}^x (xv) \cdot q_t^v (x) q_{\{s,r\}}^x (yv) \cdot q_r^y (v) q_{\{t,s\}}^v (xy),
	$$
	and hence,
	$$
	q_t^v (x) q_s^x (y) q_r^y (v) = q_s^v (y) q_t^y (x) q_r^x (v).
	$$
	The application of the reasoning used in the proof of Theorem~\ref{T-Q1f-P} and the obtained relations allows one to verify the correctness of formula~\eqref{Q1->P1} and the fact that for all $t,s \in \mathbb{Z}^d$ and $x \in X^t$, $y \in X^s$, it holds
	$$
	P_t(x) q_s^x(y) = P_s(y) q_t^y(x).
	$$
	Using this equality and the consistency conditions~\eqref{QPi-sogl}, for all $t,s \in V \subset W$ and $x \in X^t$, $y \in X^s$, $u \in X^{V \backslash \{t,s\}}$, we obtain
	$$
	P_t (x) \cdot q_{V \backslash \{t\}}^x (yu) = P_t (x) q_s^x (y) \cdot q_{V \backslash \{t,s\}}^{xy} (u) 
	= P_s (y) q_t^y (x) \cdot q_{V \backslash \{t,s\}}^{xy} (u) = P_s (y) \cdot q_{V \backslash \{s\}}^y (xu).
	$$
	Therefore, the values of $P_V$ do not depend on the choice of $t \in V$, $V \in W$, and hence, formula~\eqref{QPi->P} is also correct.
	
	It is easy to see that for all $V \in W$, the function $P_V$ defined by~\eqref{QPi->P} is a probability distribution on $X^V$. Let us show that the elements of the system $\{P_V, V \in W \}$ are consistent in Kolmogorov's sense. For all $V,I \in W$, $V \cap I = {\O}$ and $x \in X^V$, with the usage of~\eqref{QPi-sogl-I-V}, we can write
	$$
	\sum \limits_{y \in X^I} P_{V \cup I}(xy) = P_t(x_t) \sum \limits_{y \in X^I} q_{(V \backslash \{t\}) \cup I}^{x_t} (x_{V \backslash \{t\}}y) 
	= P_t(x_t) q_{V \backslash \{t\}}^{x_t} (x_{V \backslash \{t\}}) = P_V(x),
	$$
	where $t \in V$. Hence, there exists a random field $P_{{\rm Q}^{\Pi}} = \{P_V, V \in W \}$. It remains to note that according to~\eqref{QPi->P}, one has $Q^{\Pi} (P_{{\rm Q}^{\Pi}}) = {\rm Q}^{\Pi}$, and by Theorem~\ref{Th-QPPalm->P}, $P_{{\rm Q}^{\Pi}}$ is uniquely determined by ${\rm Q}^{\Pi}$.
\end{proof}

Obtained  one-to-one correspondence between a random field $P$ and a Palm specification ${\rm Q}^{\Pi}$ allows one to identify the random field $P$ with its Palm distributions and write $P = \{Q_V^z, z \in X^t, t \in \mathbb{Z}^d \backslash V, V \in W \}$. This approach, in a certain sense, is dual to the one using the system ${\rm Q}_1^f$. Indeed, the system ${\rm Q}_1^f$ consists of probability distributions indexed by one-point subsets of the lattice and parameterized by boundary conditions in finite subsets, while for the elements of ${\rm Q}^{\Pi}$, on the contrary, one-point sets of the lattice are supports for the boundary conditions and finite subsets stand for the indexes. Palm systems can be convenient for studying lattice models of point random processes.\\

The relation between $f$-- and Palm specifications is revealed in the following statements.

\begin{theorem}
	\label{Th-QPi->Qf}
	A set ${\rm Q}^{\Pi} = \left\{ q_V^z, z \in X^t, t \in \mathbb{Z}^d, V \in W(\mathbb{Z}^d \backslash \{t\}) \right\}$ of positive finite-dimensional probability distributions parameterized by boundary conditions at a point is a subsystem of an $f$--specification ${\rm Q}^f = \{q_V^z, z \in \widehat{X^f_V}, V \in W\}$ if and only if ${\rm Q}^{\Pi}$ is a Palm specification. The specification ${\rm Q}^f$ is uniquely determined by ${\rm Q}^{\Pi}$.
\end{theorem}
\begin{proof}
	Let us start with the necessity. According to~\eqref{Qf-fraction}, the elements of an $f$--specification ${\rm Q}^f = \{q_V^z, z \in \widehat{X^f_V}, V \in W\}$ satisfy the following relations
	$$
	\frac{q_{t \cup V}^y (xu)}{q_{s \cup V}^x (yu)} = \frac{q_t^y (x)}{q_s^x (y)}
	$$
	for all $V \in W$, $t,s \in \mathbb{Z}^d \backslash V$ and $x \in X^t$, $y \in X^s$, $u \in X^V$.
	These are the consistency conditions~\eqref{QPi-sogl-I-V} for the elements of Palm specification ${\rm Q}^{\Pi}$. The consistency conditions~\eqref{QPi-sogl} for the elements of ${\rm Q}^{\Pi}$ are obtained by taking the sum of both sides of~\eqref{Qf-sogl} over all $y \in X^I$ (see Remark~\ref{remark-Q-sogl-Kolm}). Hence, a subsystem of an $f$--specification, which consists of the elements parameterized by boundary conditions at a point, form the Palm specification.
	
	Now, let ${\rm Q}^{\Pi} = \{Q_V^z, z \in X^t, t \in \mathbb{Z}^d \backslash V, V \in W \}$ be a Palm specification. For any $V,\Lambda \in W$, $V \cap \Lambda = {\O}$ and $z \in X^\Lambda$, put
	\begin{equation}
		\label{QPi->Qfin}
		q_V^z(x) = \frac{q_{V \cup (\Lambda \backslash \{t\})}^{z_t} (x z_{\Lambda \backslash \{t\}})}{q_{\Lambda \backslash \{t\}}^{z_t} (z_{\Lambda \backslash \{t\}})}, \qquad x \in X^V,
	\end{equation}
	where $t \in \Lambda$. Let us show that the values of $q_V^z$ do not depend on the choice of $t \in \Lambda$. Indeed, according to the consistency conditions~\eqref{QPi-sogl}, for any $t,s \in \Lambda$, we have
	$$
	q_{V \cup (\Lambda \backslash \{t\})}^{z_t} (x z_{\Lambda \backslash \{t\}}) = \frac{q_s^{z_t}(z_s) q_{V \cup (\Lambda \backslash \{s\})}^{z_s} (x z_{\Lambda \backslash \{s\}})}{q_t^{z_s} (z_t)},
	$$
	$$
	q_{\Lambda \backslash \{t\}}^{z_t} (z_{\Lambda \backslash \{t\}}) = \frac{q_s^{z_t}(z_s) q_{\Lambda \backslash \{s\}}^{z_s} (z_{\Lambda \backslash \{s\}})}{q_t^{z_s} (z_t)}.
	$$
	Then
	$$
	\frac{q_{V \cup (\Lambda \backslash \{t\})}^{z_t} (x z_{\Lambda \backslash \{t\}})}{q_{\Lambda \backslash \{t\}}^{z_t} (z_{\Lambda \backslash \{t\}})} = \frac{q_s^{z_t}(z_s) q_{V \cup (\Lambda \backslash \{s\})}^{z_s} (x z_{\Lambda \backslash \{s\}}) q_t^{z_s} (z_t)}{q_t^{z_s} (z_t) q_s^{z_t}(z_s) q_{\Lambda \backslash \{s\}}^{z_s} (z_{\Lambda \backslash \{s\}})} = \frac{q_{V \cup (\Lambda \backslash \{s\})}^{z_s} (x z_{\Lambda \backslash \{s\}})}{q_{\Lambda \backslash \{s\}}^{z_s} (z_{\Lambda \backslash \{s\}})}.
	$$
	
	Further, according to~\eqref{QPi-sogl-I-V}, for each $V \in W$ and $z \in \widehat{X_V^f}$, the function $q_V^z$ is a probability distribution on $X^V$. To complete the proof, it remains to show that the elements of the system ${\rm Q}^f = \{q_V^z, z \in \widehat{X_V^f}, V \in W\}$ satisfy the consistency conditions~\eqref{Qf-sogl}. Let $I, V, \Lambda \in W$ be pairwise disjoint non-empty sets. For any $x \in X^V$, $y \in X^I$, $z \in X^\Lambda$ and any point $t \in \Lambda$, using~\eqref{QPi->Qfin}, we can write
	$$
	q_{V \cup I}^z (xy) = \dfrac{q_{V \cup I \cup (\Lambda \backslash {\{t\}})}^{z_t}(x y z_{\Lambda \backslash {\{t\}}})}{q_{\Lambda \backslash {\{t\}}}^{z_t} (z_{\Lambda \backslash {\{t\}}})} 
	= \dfrac{q_{V \cup (\Lambda \backslash {\{t\}})}^{z_t}(x z_{\Lambda \backslash {\{t\}}})}{q_{\Lambda \backslash {\{t\}}}^{z_t} (z_{\Lambda \backslash {\{t\}}})} \cdot \dfrac{q_{V \cup I \cup (\Lambda \backslash {\{t\}})}^{z_t}(x y z_{\Lambda \backslash {\{t\}}})}{q_{V \cup (\Lambda \backslash {\{t\}})}^{z_t}(x z_{\Lambda \backslash {\{t\}}})} = q_V^z(x) q_I^{zx}(y).
	$$
\end{proof}

\begin{theorem}
	\label{P(Qfin_1)=P(QPi)=P(Qfin)}
	Let ${\rm Q}^{\Pi}$ be a Palm specification and ${\rm Q}^f$ be the one constructed by it $f$--specification. Then $P_{{\rm Q}^{\Pi}} = P_{{\rm Q}^f}$.
\end{theorem}

\section{Systems with various boundary conditions}
\label{Section-Q}

In this section, we introduce and study systems of finite-dimensional and one-point distributions parameterized by various (both finite and infinite) boundary conditions.

\subsection{System of finite-dimensional distributions with various boundary conditions}

For a given random field $P = \{P_V, V \in W\}$, consider the set $Q_P = \{Q_V^z, z \in \widehat{X_V}, V \in W\}$ of its conditional probabilities that includes, in addition to finite-conditional probabilities $Q_P^f$, the conditional probabilities with infinite boundary conditions, determined by formula~\eqref{P->QP-lim}. By the martingale convergence theorem, the limits on the right-hand side of~\eqref{P->QP-lim} exist for almost all (in measure $P$) infinite boundary conditions. All other elements of $Q_P$ with infinite boundary conditions can be set arbitrary. Thus, for $P$, there are various systems $Q_P$, all of which will be called the \emph{full conditional distribution of the random field $P$}.

Note that in his now-classic work~\cite{D68}, Dobrushin defined a conditional distribution of a random field as a subsystem $Q_P^D$ of the system $Q_P$ consisting only of those conditional probabilities $Q_V^z$ for which $s(z) = \mathbb{Z}^d \backslash V$, $V \in W$. It seems more natural to call the system $Q_P$, which includes conditional probabilities parameterized by any (both infinite and finite) boundary conditions, a conditional distribution of a random field. However, following the tradition, we leave the term conditional distribution for the Dobrushin-type conditional distribution considered in Section~\ref{Section-QD}. The same approach is applied to the terms specification and full specification.

Since the full conditional distribution $Q_P$ contains the subsystem $Q_P^f$, which restores the random field $P$, the system $Q_P$ is also a solution to the direct problem.

\begin{theorem}
	\label{Th-QP->P}
	Any random field is restored by its full conditional distribution.
\end{theorem}

It is not difficult to see that the elements of $Q_P$ are connected by the relations~\eqref{cond-prop-main} for all finite and $P$-a.e. infinite boundary conditions. Thus, the relation~\eqref{cond-prop-main} establishes a connection between those elements of $Q_P$ whose boundary conditions differ no more than in a finite set. The relation~\eqref{P->QP-lim}, which is valid for $P$-a.e. configurations, reflects the connection between the elements with finite and infinite boundary conditions. We consider these relations as characterizing properties of $Q_P$.

The set ${\rm Q} = \{q_V^z, z \in \widehat{X_V}, V \in W\}$ of strictly positive probability distributions $q_V^z$ on $X^V$ parameterized by boundary conditions $z$ outside $V$, $V \in W$, will be called a \emph{full specification} if its elements satisfy the following consistency conditions:\\
1. for all disjoint sets $V,I\in W$ and all configurations $x \in X^V$, $y \in X^I$, $z \in \widehat{X_{V \cup I}}$, it holds
\begin{equation}
	\label{Q-sogl}
	q_{V \cup I}^z (xy) = q_V^z(x) q_I^{zx} (y);
\end{equation}
2. for all $V \in W$ and $S \subset \mathbb{Z}^d \backslash V$,
\begin{equation}
	\label{Q-sogl-Z}
	q_V^z (x) = \lim \limits_{\Lambda \uparrow S} q_V^{z_\Lambda} (x), \qquad x \in X^V, \, z \in X^S.
\end{equation}

A full specification ${\rm Q}$ is called quasilocal if its elements are quasilocal as functions on boundary conditions. Note that for a quasilocal specification ${\rm Q}$, the convergence in~\eqref{Q-sogl-Z} is uniform with respect to the boundary condition $z \in X^S$ for all $V \in W$ and $S \subset \mathbb{Z}^d \backslash V$.

Let us consider now the inverse problem of the description of random fields for a given full specification.

\begin{theorem}
	\label{Th-Q->P}
	Any full specification ${\rm Q}$ specifies a random field $P_{\rm Q}$ such that $Q_{P_Q} = {\rm Q}$ (for $P_{\rm Q}$-a.e. boundary conditions).
\end{theorem}
\begin{proof}
	Let ${\rm Q} = \{ q_V^z, z \in \widehat{X_V}, V \in W \}$ be a full specification. Using the same reasoning as in the proof of Theorem~\ref{Th-Qfin->P}, we construct a random field $P_{\rm Q} = \{P_V, V \in W\}$ such that $Q_V^z = q_V^z$ for any $V \in W$ and any finite boundary condition $z \in \widehat{X_V^f}$.
	
	Further, for any $z \in \widehat{X_V}$ and any increasing sequence of (finite) sets $\Lambda = \{\Lambda_n\}_{n \ge 1}$ such that $\Lambda \uparrow s(z)$, we have
	$$
	\frac{P_{V \cup \Lambda_n} (xz_{\Lambda_n})}{P_{\Lambda_n} (z_{\Lambda_n})} = q_V^{z_{\Lambda_n}} (x), \qquad x \in X^V, n \ge 1.
	$$
	As $n \to \infty$, the left part of the obtained relation $P_{\rm Q}$-a.e. converges to the conditional probability $Q_V^z (x)$ of the random field $P_{\rm Q}$, while its right-hand side converges to $q_V^z(x)$ by~\eqref{Q-sogl-Z}. Hence, for $P_{\rm Q}$--a.e. $z \in \widehat{X_V}$, we have
	$$
	Q_V^z (x) = \lim \limits_{n \to \infty} \frac{P_{V \cup \Lambda_n} (xz_{\Lambda_n})}{P_{\Lambda_n} (z_{\Lambda_n})} = \lim \limits_{n \to \infty} q_V^{z_{\Lambda_n}} (x) = q_V^z (x), \qquad x \in X^V, V \in W,
	$$
	that is, $Q_{P_{\rm Q}} = {\rm Q}$ for $P_{\rm Q}$-a.e. infinite boundary conditions. According to Theorem~\ref{Th-QP->P}, the random field $P_{\rm Q}$ is uniquely determined by ${\rm Q}$.
\end{proof}

Let us return to the consideration of the conditional distribution $Q_P$ of the random field $P$.
%
Since the elements of $Q_P$ are defined for $P$-a.e. boundary conditions, the random field $P$ may have many versions of its full conditional distribution. As it was mentioned above, the elements of any version of the conditional distribution satisfy relations~\eqref{Q-sogl} and~\eqref{Q-sogl-Z} for $P$-a.e. boundary conditions. Moreover, the following statement takes place.

\begin{prop}
	\label{Prop-VersCondDistr-Spec}
	For any random field $P$, there exists a version $Q_P$ of its full conditional distribution, the elements of which satisfy the consistency conditions~\eqref{Q-sogl} for all boundary conditions.
\end{prop}
\begin{proof}
	Let $Q_P= \{Q_V^z, z \in \widehat{X_V}, V \in W\}$ be a version of the full conditional distribution of the random field $P$. Denote by $\mathscr{X}_V$ the set of such configurations $z \in X^S$, $S \subset \mathbb{Z}^d \backslash V$, for which the limit in the right-hand side of~\eqref{P->QP-lim} exists, $V \in W$, and let $\mathscr{X} = \bigcup \limits_{V \in W} \mathscr{X}_V$. In this case, $P(\mathscr{X}) = 1$.
	
	Let us verify that $z \in \mathscr{X}$ if and only if $zx \in \mathscr{X}$ for any $V \in W(\mathbb{Z}^d \backslash s(z))$ and $x \in X^V$. Indeed, let $z \in \mathscr{X}$. From the definition of $\mathscr{X}$, it follows that for any disjoint sets $I,V \in W(\mathbb{Z}^d \backslash s(z))$, the following limits exist for all $x \in X^V$ and $y \in X^I$:
	$$
	\lim \limits_{\Lambda \uparrow s(z)} \frac{P_{V \cup I \cup \Lambda} (xyz_\Lambda)}{P_\Lambda (z_\Lambda)} = Q_{V \cup I}^z (xy), \qquad \lim \limits_{\Lambda \uparrow s(z)} \frac{P_{V \cup \Lambda} (xz_\Lambda)}{P_\Lambda (z_\Lambda)} = Q_V^z (x).
	$$
	But in this case, there also exists the limit
	$$
	Q_I^{zx} (y) = \lim \limits_{\Lambda \uparrow s(z)} \dfrac{P_{V \cup I \cup \Lambda} (xyz_\Lambda)}{P_{V \cup \Lambda} (xz_\Lambda)} = \lim \limits_{\Lambda \uparrow s(z)} \dfrac{P_{V \cup I \cup \Lambda} (xyz_\Lambda)}{P_\Lambda (z_\Lambda)} \cdot \lim \limits_{\Lambda \uparrow s(z)} \dfrac{P_\Lambda (z_\Lambda)}{P_{V \cup \Lambda} (xz_\Lambda)} = \dfrac{Q_{V \cup I}^z (xy)}{Q_V^z (x)}.
	$$
	Now, let $zx \notin \mathscr{X}$ for any $V \in W(\mathbb{Z}^d \backslash s(z))$ and $x \in X^V$. If $z \in \mathscr{X}$, then from the fact proved above, it follows that $zx \in \mathscr{X}$, which leads towards a contradiction. Hence, $z \notin \mathscr{X}$.
	
	Further, for all $V \in W$ and $x \in X^V$, put $q_V^z (x) = Q_V^z(x)$ if $z \in \mathscr{X}$ and $q_V^z (x) = \vert X \vert^{-\vert V \vert}$ if $z \notin \mathscr{X}$. It is clear that the system $Q = \{q_V^z, z \in \widehat{X_V}, V \in W\}$ is a version of the full conditional distribution of the random field $P$. Let us show that the elements of $Q$ satisfy the consistency conditions~\eqref{Q-sogl} for all $z$.
	
	If $z \in \mathscr{X}$, we have
	$$
	q_{V \cup I}^z (xy) = Q_{V \cup I}^z (xy) = Q_V^z(x) Q_I^{zx} (y) = q_V^z(x) q_I^{zx} (y)
	$$
	for any $I,V \in W(\mathbb{Z}^d \backslash s(z))$, $I \cap V = {\O}$ and $x \in X^V$, $y \in X^I$. In the case $z \notin \mathscr{X}$, it also holds that $zx \notin \mathscr{X}$, and hence,
	$$
	q_{V \cup I}^z (xy) = \vert X \vert^{-\vert V \cup I \vert} = \vert X \vert^{-\vert V \vert} \cdot \vert X \vert^{-\vert I \vert} = q_V^z(x) q_I^{zx} (y).
	$$
\end{proof}

In the case $P$ has a quasilocal version of its full conditional distribution, both sets of the consistency conditions hold true for all boundary conditions. Moreover, the following result takes place, that first was proved in~\cite{DN2009} for Dobrushin-type boundary conditions only (see Proposition~\ref{Prop-QDP-quasiloc} below).

\begin{prop}
	\label{Prop-QP-quasiloc}
	If a random field $P$ has a quasilocal version $Q_P$ of its full conditional distribution, then this version is unique and forms a full specification.
\end{prop}
\begin{proof}
	Let $Q_P$ be a quasilocal version of the full conditional distribution of the random field $P$. Then the limits in the right-hand side of~\eqref{P->QP-lim} exist for all boundary conditions, that is, the elements of $Q_P$ satisfy the consistency conditions~\eqref{Q-sogl-Z}. The validity of the consistency conditions~\eqref{Q-sogl} directly follows from the definition of conditional probabilities.
	
	Let us verify that the quasilocal version is unique. Assume the opposite. Let $Q_P = \{ Q_V^z, z \in \widehat{X_V}, V \in W \}$ and $F_P = \{ F_V^z, z \in \widehat{X_V}, V \in W \}$ be two quasilocal versions of the full conditional distribution of the random field $P$. Then for all $V \in W$ and $x \in X^V$, the function
	$$
	f_S(z) = Q_V^z (x) - F_V^z (x), \qquad z \in X^S, S \subset \mathbb{Z}^d \backslash V,
	$$
	equals to zero for $P$-a.e. $z \in X^S$, i.e., 
	$$
	P(X^S_0) = P(\{z \in X^S: \, f_S(z) \neq 0\}) = 0.
	$$
	Further, from the quasilocality of $Q_P$ and $F_P$, it follows the quasilocality of $f_S$, and thus, for any $\varepsilon > 0$, there exists $\Lambda_0 \in W(S)$ such that for all $\Lambda \supset \Lambda_0$, $\Lambda \in W(S)$, it holds $\left| f_S(z) - f_\Lambda(z_\Lambda) \right| < \varepsilon$. Hence, if $z \in X^S_0$, there exists $\Lambda \in W(S)$ large enough such that $f_\Lambda(z_\Lambda) \neq 0$. But in this case, by positivity of $P$,
	$$
	P(X^S_0) \ge P(X^\Lambda_0) \ge P(z_\Lambda) > 0,
	$$
	and we come to the contradiction.
\end{proof}

It is easy to see that the subsystem of a full specification ${\rm Q}$, which consists only of those elements which are parameterized by finite boundary conditions, forms an $f$--specification ${\rm Q}^f$. On the other hand, an $f$--specification, generally speaking, does not restore the full specification for which it is a subsystem. However, the following statement holds true.

\begin{prop}
	\label{Prop-Qf->Q}
	Let ${\rm Q}^f$ be an $f$--specification. Then there exists the system $Q^* = \{q_V^z, z \in \widehat{X_V}, V \in W\} $ of positive probability distributions parameterized by various boundary conditions, all the elements of which satisfy the consistency conditions~\eqref{Q-sogl}. In this case, ${\rm Q}^f \subset {\rm Q}^*$.
\end{prop}
\begin{proof}
	According to Theorem~\ref{Th-QPfin->P}, for a given $f$--specification ${\rm Q}^f$, there exists a unique random field $P$ such that $Q^f_P = {\rm Q}^f$. By Proposition~\ref{Prop-VersCondDistr-Spec}, this random field $P$ has a version ${\rm Q}_P$ of its full conditional distribution which elements satisfy the consistency conditions~\eqref{Q-sogl} for all boundary conditions. It remains to put ${\rm Q}^* = {\rm Q}_P$.
\end{proof}

Note also, that for the elements of full specification ${\rm Q}$, Remarks 1--3 are also stay true.

\subsection{System of one-point distributions parameterized by various boundary conditions}

The system $Q_1(P) = \{Q_t^z, z \in \widehat{X_t}, t \in \mathbb{Z}^d \}$ of one-point conditional probabilities of a random field $P$, defined by formulas~\eqref{P->QP} and~\eqref{P->QP-lim} for $V = \{t\}$, $t \in \mathbb{Z}^d$, will be called a \emph{full one-point conditional distribution of the random field $P$}.

It is clear that $Q_1^f(P) \subset Q_1(P)$, and hence, full one-point conditional distribution of a random field is a solution to the direct problem.

\begin{theorem}
	\label{Th-Q1P->QP->P}
	Any random field is restored by its full one-point conditional distribution.
\end{theorem}

As the main characterizing properties of one-point conditional probabilities we consider the property~\eqref{cond1-prop-main-2x2} as well as the relation~\eqref{P->QP-lim} (for $V = \{t\}$, $t \in \mathbb{Z}^d$) establishing a connection between finite and infinite boundary conditions.

A system ${\rm Q}_1 = \{q_t^z, z \in \widehat{X_t}, t \in \mathbb{Z}^d \}$ of strictly positive one-point probability distributions parameterized by various boundary conditions will be called a \emph{full 1--specification} if its elements satisfy the following consistency conditions:\\
1. for all $t,s \in \mathbb{Z}^d$ and $x \in X^t$, $y \in X^s$, $z \in \widehat{X_{\{t,s\}}}$, it holds
\begin{equation}
	\label{Q1-sogl-2x2}
	q_t^z (x) q_s^{zx} (y) = q_s^z (y) q_t^{zy} (x);
\end{equation}
2. for any $t \in \mathbb{Z}^d$ and $S \subset \mathbb{Z}^d \backslash \{t\}$,
\begin{equation}
	\label{Q1-sogl-Z}
	q_t^z (x) = \lim \limits_{\Lambda \uparrow S} q_t^{z_\Lambda} (x), \qquad x \in X^t, \, z \in X^S.
\end{equation}

A full 1--specification ${\rm Q}_1$ is called quasilocal if its elements are quasilocal as functions on boundary conditions. For a quasilocal full 1--specification ${\rm Q}_1$, the convergence in~\eqref{Q1-sogl-Z} is uniform in $z \in X^S$ for all $t \in \mathbb{Z}^d$ and $S \subset \mathbb{Z}^d \backslash \{t\}$.

Let us consider the inverse problem of the description of random fields for a full 1--specification.

\begin{theorem}
	\label{T-Q1-P}
	Any full 1--specification ${\rm Q}_1$ specifies a random field $P_{{\rm Q}_1}$ such that $Q_1 (P_{{\rm Q}_1}) = {\rm Q}_1$ (for $P_{{\rm Q}_1}$-a.e. boundary conditions).
\end{theorem}

The proof of this result is similar to the proof of Theorem~\ref{Th-Q->P} (using Theorem~\ref{T-Q1f-P}) and, therefore, will be omitted.

Note that the elements of full 1--specification ${\rm Q}_1$ satisfy the following conditions: for any $t,s \in \mathbb{Z}^d$ and $x,u \in X^t$, $y,v \in X^s$, $z \in \widehat{X^0_{\{t,s\}}}$, the following relations hold
\begin{equation}
	\label{Q1-sogl}
	q_t^{zy} (x) q_s^{zx} (v) q_t^{zv} (u) q_s^{zu} (y) = q_t^{zy} (u) q_s^{zu} (v) q_t^{zv} (x) q_s^{zx} (y).
\end{equation}

It is not difficult to see that the one-point subsystem of a full specification forms a full 1--specification. Moreover, the following statement is true, which can be shown analogously to the proof of Theorem~\ref{Th-Q1f->Qf}.

\begin{theorem}
	\label{Th-Q1->Q}
	A set ${\rm Q}_1 = \{q_t^z, z \in \widehat{X_t}, t \in \mathbb{Z}^d\}$ of positive one-point probability distributions parameterized by various boundary conditions is a one-point subsystem of a full specification ${\rm Q} = \{q_V^z, z \in \widehat{X_V}, V \in W\}$ if and only if ${\rm Q}_1$ is a full 1--specification. The specification ${\rm Q}$ is uniquely determined by ${\rm Q}_1$.
\end{theorem}

As in the case of specifications with finite boundary conditions, the full specification ${\rm Q}$ can be constructed by the elements of the full 1--specification ${\rm Q}_1$ using either formula~\eqref{Q1->Q-1} or formula~\eqref{Q1->Q-2}. The full specification ${\rm Q}$ constructed from ${\rm Q}_1$ inherits such properties of ${\rm Q}_1$ as being quasilocal or Markovian. Moreover, the following fact takes place.

\begin{theorem}
	Let ${\rm Q}_1$ be a full 1--specification and ${\rm Q}$ be the constructed by it full specification. Then $P_{{\rm Q}_1} = P_{\rm Q}$.
\end{theorem}

Concluding this section, we note several properties of the full one-point conditional distribution of a random field. These results directly follow from the similar statements for a full conditional distribution or can be verified independently using similar reasoning.


\begin{prop}
	\label{Prop-VersCondDistr1-1Spec}
	For any random field $P$, there exists a version $Q_1(P)$ of its full one-point conditional distribution, the elements of which satisfy the relations~\eqref{cond1-prop-main-2x2} for all boundary conditions.
\end{prop}

\begin{prop}
	\label{Prop-Q1P-quasiloc}
	If a random field $P$ has a quasilocal version $Q_1(P)$ of its full one-point conditional distribution, this version is unique and forms a full 1--specification.
\end{prop}

\section{Systems of Dobrushin-type conditional distributions}
\label{Section-QD}

Among the subsystems of the full conditional distribution of a random field, the system introduced by Dobrushin in~\cite{D68} occupies a special place. Interest in this system is caused, first of all, by the problems of mathematical statistical physics.

Dobrushin was the first to consider the problem of the description of a random field by conditional probabilities. Further, Dachian and Nahapetian in the series of works~\cite{DN1998, DN2001, DN2004} showed that Dobrushin's theory can be equivalently formulated in terms of the system of consistent one-point distributions parameterized by boundary conditions.

In this section, we formulate the main results of the mentioned works from the point of view developed in the present paper.

\subsection{System of finite-dimensional distributions with infinite boundary conditions}

For a random field $P$, considered by Dobrushin~\cite{D68} the system $Q^D_P = \{Q_V^z, z \in X^{\mathbb{Z}^d \backslash V}, V \in W\}$ of conditional probabilities on $X^V$ parameterized by infinite boundary conditions defined everywhere outside $V$, $V \in W$, will be called \emph{infinite conditional distribution of the random field $P$}, or just \emph{conditional distribution of $P$}.

Dobrushin's system is not a solution to the direct problem, since different random fields can have the same infinite conditional distribution (see, for example,~\cite{DN2019}). Nevertheless, one can single out a class of random fields that can be restored by their conditional distribution (see Theorem 2 in~\cite{D68}).

\begin{theorem}
	Let a random field $P$ be such that its conditional distribution $Q^D_P$ is quasilocal and satisfy the following condition
	\begin{equation}
		\label{UniqCond}
		\sum \limits_{s \in \mathbb{Z}^d \backslash \{t\}} {\rho _{s,t}} \le a < 1,
	\end{equation}
	where
	$$
	\rho_{s,t}  = \sup \limits_{z,y \in X^{\mathbb{Z}^d \backslash \{ t \}}: z_{\mathbb{Z}^d \backslash \{ t,s \}} = y_{\mathbb{Z}^d \backslash \{ t,s \}}} \frac{1}{2} \sum \limits_{x \in X^t} \left| Q_t^z (x) - Q_t^y (x) \right|.
	$$
	Then the system $Q^D_P$ restores $P$.
\end{theorem}

The main property of conditional probabilities considered by us cannot be written for Dobrushin-type conditional probabilities directly in the form~\eqref{cond-prop-main}. However, as it is not difficult to see, for the elements of $Q^D_P$, relations~\eqref{cond-prop-main} can be written as follows
$$
Q_{V \cup I}^z (xy) = Q_I^{zx} (y) \sum \limits_{\beta \in X^I} Q_{V \cup I}^z (x\beta),
$$
where $V,I \in W$, $V \cap I = {\O}$ and $x \in X^V$, $y \in X^I$, $z \in X^{\mathbb{Z}^d \backslash (V \cup I)}$. It is this relation that Dobrushin considered as the characterizing property of conditional probabilities with infinite boundary conditions.

The set ${\rm Q}^D = \{q_V^z, z \in X^{\mathbb{Z}^d \backslash V}, V \in W\}$ of strictly positive probability distributions parameterized by infinite boundary conditions will be called a \emph{specification in Dobrushin's sense}, or just \emph{specification} if its elements satisfy the following consistency conditions: for all disjoint sets $V,I \in W$ and all configurations $x \in X^V$, $y \in X^I$, $z \in X^{\mathbb{Z}^d \backslash (V \cup I)}$, it holds
\begin{equation}
	\label{QD-sogl}
	q_{V \cup I}^z (xy) = q_I^{zx} (y) \sum \limits_{\beta \in X^I} q_{V \cup I}^z (x\beta).
\end{equation}

Dobrushin presented conditions under which a specification ${\rm Q}^D$ defines a random field (see Theorem 1 in~\cite{D68}). In this case, there may exist several random fields whose conditional distribution a.e. coincides with ${\rm Q}^D$. However, the conditions on the elements of ${\rm Q}^D$ under which it specifies (uniquely determines) compatible with random field are known (see, for example, Theorem 2 in~\cite{D68}).

\begin{theorem}
	Let ${\rm Q}^D$ be a quasilocal specification. Then there exists a random field $P_{{\rm Q}^D}$ such that $Q^D (P_{{\rm Q}^D}) = {\rm Q}^D$ ($P_{{\rm Q}^D}$-a.e.). If, in addition, condition~\eqref{UniqCond} is satisfied, then the random field $P_{{\rm Q}^D}$ is unique.
\end{theorem}

Note that in the theory of Gibbs random fields, the inverse problem for Dobrushin's specification ${\rm Q}^D$ is usually formulated in terms of DLR--equations. Namely, for a given (Gibbs) specification $Q^D = \{q_V^z, z \in X^{\mathbb{Z}^d \backslash V},\linebreak V \in W\}$, one considers the question of the existence as well as the uniqueness of a random field $P=\{P_V, V \in W\}$ satisfying the following equations
$$
P_V(x) = \int \limits_{z \in X^{\mathbb{Z}^d \backslash V}} q_V^z(x) P_{\mathbb{Z}^d \backslash V}(z)
$$
for all $x \in X^V$ and $V \in W$.

\begin{remark}
	The consistency conditions~\eqref{QD-sogl} can be written in the following equivalent form: for all disjoint sets $V,I \in W$ and all configurations $x, u \in X^V$, $y \in X^I$, $z \in X^{\mathbb{Z}^d \backslash (V \cup I)}$, it holds
	$$
	q_{V \cup I}^z (xy) q_V^{zy}(u) = q_{V \cup I}^z (uy) q_V^{zy} (x).
	$$
	These relations are true for the elements of $f$ and full specification as well (see Remark~\ref{R2}).
\end{remark}

\begin{remark}
	The elements of a specification ${\rm Q}^D$ are connected with each other by the following analogues of relations~\eqref{qz-qzy}: for all $I,V \in W$, $V \cap I = {\O}$, and $z \in X^{\mathbb{Z}^d \backslash (V \cup I)}$, it holds
	$$
	\sum \limits_{\beta \in X^I} q_{V \cup I}^z (x \beta) = \frac{q_V^{zy}(x)}{q_I^{zx}(y)} \left(\sum \limits_{\alpha \in X^V} \frac{q_V^{zy}(\alpha)}{q_I^{z \alpha}(y)} \right)^{-1}, \qquad x \in X^V,
	$$
	where $y \in X^I$.
\end{remark}

Any full specification ${\rm Q}$ contains a subsystem ${\rm Q}^D$ which is a specification in the sense of Dobrushin. However, not every specification ${\rm Q}^D$ defines a full specification ${\rm Q}$ for which it is a subsystem. Moreover, if ${\rm Q}^D$ defines some full specification ${\rm Q}$ and some random field $P_{{\rm Q}^D}$, it may turn out that $P_{\rm Q} \neq P_{{\rm Q}^D}$.

Let us mention some properties of conditional distribution $Q^D_P$ of a random filed $P$. It is clear, that $P$ may have several versions of its infinite conditional distribution $Q^D_P$. Also, the following statements take place (see, for example, Theorem in~\cite{Gold} and Proposition 3.3 in~\cite{DN2009}, and compare with Propositions~\ref{Prop-VersCondDistr1-1Spec} and~\ref{Prop-Q1P-quasiloc} of the present paper).


\begin{prop}
	\label{Prop-QDP-everywhere}
	For a random field $P$, there exists a version $Q^D_P$ of its conditional distribution which forms a specification.
\end{prop}

\begin{prop}
	\label{Prop-QDP-quasiloc}
	If a random field $P$ has a quasilocal version $Q^D_P$ of its conditional distribution, then this version is unique and forms a specification.
\end{prop}

Dobrushin~\cite{D68} defined a Markov random field $P$ (with respect to a neighborhood system $\partial = \{ \partial t, t \in \mathbb{Z}^d \}$ in $\mathbb{Z}^d$) as one for which the elements of $Q_P^D$ satisfy the following conditions: for all $V \in W$ and $P$-a.e. $z \in X^{\mathbb{Z}^d \backslash V}$, it holds
\begin{equation}
	\label{Markov-Dobr}
	Q_V^z (x) = Q_V^{z_{\partial V}} (x), \qquad x \in X^V.
\end{equation}
Let us show that this definition is equivalent to the one given in Section~\ref{Section-Qf}.

\begin{prop}
	\label{Prop-Markov-ekv}
	A random field $P$ is Markovian if and only if the elements of its conditional distribution $Q_P^D$ satisfy the conditions~\eqref{Markov-Dobr}.
\end{prop}
\begin{proof}
	Let $P$ be a Markov field. Then for any disjoint sets $V,\Lambda \in W$ such that $\partial V \subset \Lambda$, and any configuration $z \in X^{\mathbb{Z}^d \backslash V}$, we have
	$$
	Q_V^{z_\Lambda}(x) = Q_V^{z_{\partial V}}(x), \qquad x \in X^V.
	$$
	Passing in this relation to the limit as $\Lambda \uparrow \mathbb{Z}^d \backslash V$, we obtain
	$$
	Q_V^z (x) = Q_V^{z_{\partial V}} (x) \; (P\text{-a.e.}), \qquad x \in X^V,
	$$
	and hence, the elements of conditional distribution $Q_P^D$ of the random field $P$ satisfy relations~\eqref{Markov-Dobr}.
	
	Now, let a random field $P$ be such that the elements of its conditional distribution $Q_P^D$ satisfy~\eqref{Markov-Dobr}. By Sullivan's inequality~\eqref{Q-ineq}, for any $V, \Lambda \in W$, $V \cap \Lambda = {\O}$ and $x \in X^V$, $z \in X^\Lambda$, we have
	$$
	\inf \limits_{y \in X^{\mathbb{Z}^d \backslash V}: y_\Lambda = z} Q_V^y (x) \le Q_V^z (x) \le \sup \limits_{y \in X^{\mathbb{Z}^d \backslash V}: y_\Lambda = z} Q_V^y (x).
	$$
	Thus, if $\Lambda$ is such that $\partial V \subset \Lambda$, we obtain
	$$
	\sup \limits_{y \in X^{\mathbb{Z}^d \backslash V}: y_\Lambda = z} Q_V^y (x) = \sup \limits_{y \in X^{\mathbb{Z}^d \backslash V}: y_\Lambda = z} Q_V^{y_{\partial V}} (x) = Q_V^{z_{\partial V}} (x),
	$$
	$$
	\inf \limits_{y \in X^{\mathbb{Z}^d \backslash V}: y_\Lambda = z} Q_V^y (x) = \inf \limits_{y \in X^{\mathbb{Z}^d \backslash V}: y_\Lambda = z} Q_V^{y_{\partial V}} (x) = Q_V^{z_{\partial V}} (x),
	$$
	and hence, $Q_V^z (x) = Q_V^{z_{\partial V}} (x)$. Therefore, $P$ is a Markov random field.
\end{proof}

\subsection{System of Dobrushin-type one-point distributions}

Considered in~\cite{DN1998, DN2001, DN2004} system $Q^D_1 (P) = \{Q_t^z, z \in X^{\mathbb{Z}^d \backslash \{t\}}, t \in \mathbb{Z}^d\}$ of one-point conditional distributions with infinite boundary conditions will be called the \emph{Dobrushin-type one-point conditional distribution of the random field $P$}, or, in short, the \emph{one-point conditional distribution of $P$}.

The distribution $Q^D_1 (P)$, generally speaking, does not restor the random field $P$ (see the corresponding remarks for the system $Q^D_P$).

\begin{theorem}
	Let a random field $P$ be such that its one-point conditional distribution $Q^D_1(P)$ is quasilocal and satisfy the condition~\eqref{UniqCond}. Then the system $Q^D_1(P)$ restores~$P$.
\end{theorem}

The main characterizing property~\eqref{cond1-prop-main-2x2} of one-point conditional probabilities cannot be written down directly for the elements of $Q^D_1(P)$. However, as it is not difficult to verify, for the one-point conditional probabilities of Dobrushin's type, it holds
$$
Q_t^{zy} (x) Q_s^{zx} (v) Q_t^{zv} (u) Q_s^{zu} (y) = Q_t^{zy} (u) Q_s^{zu} (v) Q_t^{zv} (x) Q_s^{zx} (y)
$$
for all $t,s \in \mathbb{Z}^d$, $x,u \in X^t$, $y,v \in X^s$ and $P$-a.e. boundary conditions $z \in X^{\mathbb{Z}^d \backslash \{t,s\}}$. It is this relation that was singled out in works~\cite{DN2001, DN2004} as the determining one for the systems of one-point distributions parameterized by Dobrushin's type boundary conditions.

A set ${\rm Q}_1^D = \{q_t^z, z \in X^{\mathbb{Z}^d \backslash \{t\}}, t \in \mathbb{Z}^d\}$ of strictly positive one-point probability distributions parameterized by infinite boundary conditions will be called a \emph{1--specification} (\emph{in the sense of Dobrushin}) if its elements satisfy the following consistency conditions: for all $t,s \in \mathbb{Z}^d$ and $x,u \in X^t$, $y,v \in X^s$, $z \in X^{\mathbb{Z}^d \backslash \{t,s\}}$, it holds
\begin{equation}
	\label{Q1D-sogl}
	q_t^{zy} (x) q_s^{zx} (v) q_t^{zv} (u) q_s^{zu} (y) = q_t^{zy} (u) q_s^{zu} (v) q_t^{zv} (x) q_s^{zx} (y).
\end{equation}

The inverse problem for a 1--specification may not have a solution. However, the following statement holds true (see Theorem 4.2 in~\cite{DN2001} and Theorem 21 in~\cite{DN2004}).

\begin{theorem}
	Let ${\rm Q}_1^D$ be a quasilocal 1--specification. Then there exists a random field $P_{{\rm Q}_1^D}$ such that $Q^D_1(P_{{\rm Q}_1^D}) = {\rm Q}^D_1$ ($P$-a.e.). If, in addition, the condition~\eqref{UniqCond} is satisfied, then the random field $P_{{\rm Q}_1^D}$ is unique.
\end{theorem}

The connection between Dobrushin's type 1--specification and specification is given in the following statements (see Theorem 19 in~\cite{DN2004}).

\begin{prop}
	Let ${\rm Q}^D_1$ be a 1--specification. Then there exists a unique specification ${\rm Q}^D$ such that ${\rm Q}^D_1 \subset {\rm Q}^D$.
\end{prop}

The proof of this result was obtained in~\cite{DN2004}, where the construction of the elements of specification ${\rm Q}^D$ by the elements of 1--specification ${\rm Q}^D_1$ was carried out according to formula~\eqref{Q1->Q-2}. From this formula, in particular, it follows that the constructed specification ${\rm Q}^D$ inherits such properties of ${\rm Q}^D_1$ as quasilocality and Markovness. Moreover, the following fact holds true.

\begin{theorem}
	\label{P(Q_1)=P(Q)}
	Let ${\rm Q}_1^D$ be a 1--specification and let ${\rm Q}^D$ be the constructed by it specification. Then the set of random fields compatible with ${\rm Q}_1^D$ coincides with the set of random fields compatible with ${\rm Q}^D$.
\end{theorem}

At the same time, 1--specification, generally speaking, does not define a full specification ${\rm Q}_1$ such that ${\rm Q}_1^D \subset {\rm Q}_1$.

Returning to the consideration of the one-point conditional distribution of a random field, we note that the random field $P$ may have many versions of it, while the following statements are true (see Propositions~\ref{Prop-VersCondDistr-Spec} and~\ref{Prop-QP-quasiloc} in this paper and Proposition 3.3 in~\cite{DN2009}).


\begin{prop}
	For a random field $P$, there exists a version $Q^D_1(P)$ of its one-point conditional distribution which forms a 1--specification.
\end{prop}

\begin{prop}
	If a random field $P$ has a quasilocal version $Q^D_1(P)$ of its one-point conditional distribution, then this version is unique and forms a 1--specification.
\end{prop}

It should be noted that the system of one-point Dobrushin type conditional distributions is a fundamental object, in terms of which the foundations of the general theory of Gibbs random fields were laid (see~\cite{DN2009} and~\cite{DN2019}).

Note also, that if one considers a finite-volume 1--specification $Q_1^\Lambda = \{q_t^z, z \in X^{\Lambda \backslash \{t\}}, t \in \Lambda\}$ as a system of probability distributions consistent in the sense of~\eqref{Q1D-sogl}, then it specifies compatible with it finite random field $P_\Lambda$, $\Lambda \in W$. For details, see~\cite{KhN22}.

\section*{Concluding remarks}

We considered various ways of the description of random fields by systems of consistent finite-dimensional distributions parameterized by boundary conditions. We presented (in the majority of cases, necessary and sufficient) conditions on the elements of these systems to coincide with the corresponding conditional probabilities of a random field.

The preference for one or another system, of course, depends on the task at hand. First of all, we note that the Dobrushin-type systems $Q^D$ and $Q^D_1$ are widely used in mathematical problems of statistical physics. Despite the fact that these systems do not specify a random field, perhaps, their main role is revealed in the theory of phase transitions. Namely, in Dobrushin's theory of description of Gibbs random fields, the  non-uniqueness of the solution to the inverse problem is interpreted as the presence of a phase transition in the model under consideration (see the fundamental work~\cite{DUniq}).

The system $Q^f$ of distributions parameterized by finite boundary conditions uniquely determines a random field and, therefore, can be useful in many theoretical considerations. It is especially convenient that it turns out to be sufficient to have a one-point system $Q^f_1$ of such distributions. First of all, we note that this system seems to be the most natural for application in the theory of Markov random fields. In particular, it can be used to describe Gaussian Markov random fields that are ubiquitous in various applications (corresponding paper is being prepared; see also~\cite{Kunc}). In addition, many properties of a random field are expressed in terms of conditions on $Q^f_1$. Note, for example, the mixing conditions~\cite{DN2011}, the constructive uniqueness criterion~\cite{DSh}, and the fact that the probabilistic definition of a Gibbs random field (without using the notion of potential) was given in~\cite{DN2009}  precisely in terms of the system $Q^f_1$. As regards the system $Q^{\Pi}$ of Palm-type distributions, it seems to be useful in studying discrete models of point processes.

The system $Q$ ($Q_1$), being the most general, is primarily of theoretical interest. As already noted, it is natural to call the system $Q_P$ the conditional distribution of the random field $P$. On the other hand, despite its generality, it admits a convenient representation: its elements can be represented in the Gibbsian form (this issue will be considered in a separate publication).

As for the main property of conditional probabilities used as the consistency conditions (see the relations~\eqref{Qf-sogl} and~\eqref{Q-sogl}, which differ only in restrictions on the supports of the boundary conditions), we note that this relation was used by Renyi~\cite{R} as the third axiom in constructing his axiomatic of probability theory by means of conditional probabilities. We also note that in~\cite{AN}, it was shown how, based on the relation~\eqref{Q-sogl}, to obtain the consistency conditions for other systems of probability distributions considered in the paper.

For the sake of simplicity, in this paper, we considered systems of strictly positive probability distributions and positive random fields only. However, the similar results can be obtained under less restrictive positivity conditions applying the approach introduced in~\cite{DN2004} (see also~\cite{KhN22}). Our results carry over in a natural way to the case of infinite (both countable and continuous) measurable spaces $X$ (under a suitable integrability condition).

Finally, let us note that both direct and inverse problems of the description of random fields can be solved in terms of systems of correlation functions (see~\cite{KhN21-Lob} and the references therein). Also, this problem can be considered from the algebraic point of view as a problem of consistency of an appropriate infinite system of linear equations (see~\cite{KhN19}). \\

\textbf{Acknowledgements.} The author is grateful to Prof. Boris S. Nahapetian for fruitful discussions and helpful suggestions.

\end{document}